# Central limit theorems for double Poisson integrals

GIOVANNI PECCATI[1] and MURAD S. TAQQU[2]

[1]*Laboratoire de Statistique Théorique et Appliquée, Université Paris VI, France.*
*E-mail: giovanni.peccati@gmail.com*
[2]*Boston University, Department of Mathematics, 111 Cummington Street, Boston, MA 02215, USA. E-mail: murad@math.bu.edu*

Motivated by second order asymptotic results, we characterize the convergence in law of double integrals, with respect to Poisson random measures, toward a standard Gaussian distribution. Our conditions are expressed in terms of contractions of the kernels. To prove our main results, we use the theory of stable convergence of generalized stochastic integrals developed by Peccati and Taqqu. One of the advantages of our approach is that the conditions are expressed directly in terms of the kernel appearing in the multiple integral and do not make any explicit use of asymptotic dependence properties such as mixing. We illustrate our techniques by an application involving linear and quadratic functionals of generalized Ornstein–Uhlenbeck processes, as well as examples concerning random hazard rates.

## 1. Introduction

Let $\widehat{N} = \{\widehat{N}(B) : B \in \mathcal{Z}\}$ be a centered Poisson measure over some Borel space $(Z, \mathcal{Z})$ and let $\{g_n : n \geq 1\}$ and $\{f_n : n \geq 1\}$ be, respectively, a sequence of real-valued functions over $Z$ and a sequence of real-valued symmetric functions over $Z^2$, vanishing on the diagonal set $\{(a, b) \in Z^2 : a = b\}$. The aim of this paper is to characterize the convergence in distribution, toward a bivariate Gaussian law, of sequences of random vectors of the type

$$\{I_1(g_n), I_2(f_n) : n \geq 1\}, \tag{1}$$

where

$$I_1(g_n) = \int_Z g_n(a) \widehat{N}(\mathrm{d}a), \qquad n \geq 1, \tag{2}$$

are single Wiener–Itô stochastic integrals with respect to $\widehat{N}$ and

$$I_2(f_n) = \int_Z \int_Z f_n(a,b)\widehat{N}(\mathrm{d}a)\widehat{N}(\mathrm{d}b), \qquad n \geq 1, \tag{3}$$

is a double Wiener–Itô integral with respect to $\widehat{N}$ (see [24], as well as the discussion contained in Section 2, for precise definitions). One of the main motivations for the study of random objects such as (2) and (3) is that the (Hilbert) space of random variables $Z$ having the representation

$$Z = c + I_1(g) + I_2(f), \tag{4}$$

for some real constant $c$ and some (uniquely defined) deterministic kernels $g$ and $f$, coincides with the $L^2$-space generated by random variables of the type $\mathbf{P}(\widehat{N}(A_1), \ldots, \widehat{N}(A_m))$, $m \geq 1$, where $\mathbf{P}(\cdot)$ is a polynomial of degree $\leq 2$ in $m$ variables (see, e.g., Ogura [14], or Rota and Wallstrom [21], for a proof of this fact and for a representation of multiple Poisson integrals in terms of Charlier polynomials). It follows that a characterization of the asymptotic normality of a sequence such as (1) is crucial in the study of the asymptotic behavior of random variables that are linear or quadratic transformations of the measure $\widehat{N}$. For instance, we will see below that one can effectively study *linear and quadratic functionals* of stochastic processes constructed from $\widehat{N}$ by first resorting to their representation in the form (4).

The main point in the characterization of the asymptotic normality of (1) is to obtain criteria for the weak convergence (toward a Gaussian law) of the sequence $\{I_2(f_n) : n \geq 1\}$. To do this, we shall make extensive use of the results, involving 'generalized adapted stochastic integrals,' developed in [17] and [19]. These results are based on a decoupling technique, known as the *principle of conditioning* (see, e.g., Xue [27]). In particular, the sufficient conditions for the asymptotic normality of $\{I_2(f_n) : n \geq 1\}$ will be expressed in terms of the *contraction kernels* associated with the functions $f_n$. These contraction kernels, whose definition is given in Section 2 below, can be easily computed from the analytic expression for $f_n$.

We will also show that if the variances of $I_2(f_n)$ converge to 1 and if some integrability assumptions are satisfied (e.g., if the sequence $\{I_2(f_n)^4 : n \geq 1\}$ is uniformly integrable), our conditions are necessary and sufficient, and also that they are equivalent to the convergence of the fourth moments $\mathbb{E}[I_2(f_n)^4]$ toward 3. Note that the number 3 is simply the fourth moment of a standard Gaussian random variable: in our opinion, this fact is rather striking, as it shows that central limit theorems (CLTs) involving the integrals $I_2(f_n)$, $n \geq 1$, are essentially determined by the convergence of quantities associated with their first two even moments. This implication should be regarded as a drastic simplification of the method of *moments and cumulants*, customarily adopted to prove CLTs for polynomial forms in random variables with finite moments of any order. See the two surveys by Surgailis, [24] and [25], for more details on this point.

Our findings extend to the framework of Poisson measures results previously established in [13, 18] and [20] for sequences of multiple Wiener–Itô integrals (of arbitrary order) with respect to general Gaussian processes. For instance, in the Gaussian case, the necessary and sufficient conditions for a CLT involving multiple integrals also require the mere convergence of second and fourth moments and can equivalently be expressed in terms of contraction kernels. Note that the results proved in [13] and [20] have already been applied in a variety of frameworks: see, for example, [9] for applications to self-intersection local times of fractional Brownian motion and [12] for an application to high-frequency CLTs on commutative groups.

As an illustration of our techniques, in Section 4 we focus on *generalized Volterra processes* associated with independently scattered random measures. These processes have the form

$$Y_t^h = \int_\mathbb{R} \int_\mathbb{R} uh(t,s)\widehat{N}(\mathrm{d}u, \mathrm{d}s), \qquad t \geq 0, \tag{5}$$

where $\widehat{N}$ is a suitable Poisson measure over $\mathbb{R} \times \mathbb{R}$ and $h$ is deterministic. The function $h$ is called the *kernel* of $Y^h$. When $h$ has the form $h(t,s) = v(t-s)$, where $v$ has support in $\mathbb{R}_+$, then $Y^h$ is called a (stationary) *moving average Lévy process*. To keep the length of this paper within bounds, we only treat the specific example of a *Ornstein–Uhlenbeck Lévy process*, which is a moving average process corresponding to the case $h(t,s) = \sqrt{2\lambda}\exp(-\lambda(t-s))\mathbf{1}_{t>s}$, $\lambda > 0$. These processes appear, for example, in Bayesian survival analysis (see James [10]), network modeling (see Cohen and Taqqu [4] or Wolpert and Taqqu [26]) and finance (see Bandorff-Nielsen and Shepard [3]). We will characterize the asymptotic behavior of the functionals

$$V_T = \frac{1}{T}\int_0^T (Y_t^h)^2 \,\mathrm{d}t, \qquad T > 0,$$

as $T \to \infty$. We do so by expressing the centered random variable $V_T - T^{-1}\mathbb{E}\int_0^T (Y_t^h)^2\,\mathrm{d}t$ as the sum of a single and a double integral with respect to $\widehat{N}$, and we use our results to prove the existence of positive constants $C(T)$ such that

$$C(T)\left[V_T - \frac{1}{T}\mathbb{E}\int_0^T (Y_t^h)^2 \,\mathrm{d}t\right] \overset{\text{law}}{\to} \mathcal{N}(0,1), \tag{6}$$

where $\mathcal{N}(0,1)$ is a standard Gaussian random variable.

This approach can be used beyond the framework of Ornstein–Uhlenbeck processes. Indeed, further applications of our main results have been developed in the two papers [5] and [15], where several CLTs are proved, involving linear and quadratic functionals of *random hazard rates* in nonparametric Bayesian survival analysis. Random hazard rates can be thought of as the mathematical representation of the 'instantaneous risk' associated with the average lifetime of a given population. In a Bayesian framework, they are modeled in terms of moving averages of the type (5). In particular, the results proved in [5] and [15] can be used in prior specification procedures. See Section 4.2 for more details on this subject.

As discussed below and in [5] and [15], one of the advantages of our techniques is that the conditions one has to check in order to prove that (6) holds can be expressed directly in terms of the kernel $h$ appearing in (5) and do not make any explicit use of (asymptotic) dependence properties of the process $Y_t^h$ (such as mixing).

Observe that the theory developed in this paper also applies to a wider class of quadratic forms associated with $Y_t^h$ such as

$$\int_0^T \int_0^T Y_s^h Y_t^h \alpha(\mathrm{d}s, \mathrm{d}t) \quad \text{and} \quad \sum_{1 \leq j, i \leq n} m_{i,j} Y_i^h Y_j^h,$$

where $\alpha$ and $m$ are, respectively, a measure and a real matrix. This kind of application, which will be the object of a separate study, should be compared with the findings contained in [1], [2] and [6].

The paper is organized as follows. In Section 2, we define multiple stochastic integrals and discuss some of their basic properties. In Section 3, we focus on sequences of single and double Poisson integrals, for which a general central limit theorem is proved. In Section 4, we describe some applications involving Ornstein–Uhlenbeck processes, random hazard rates and Bayesian survival analysis.

## 2. Definitions and preliminary results

Throughout the paper, $(Z, \mathcal{Z}, \mu)$ is a Borel measure space, with $\mu$ non-atomic, $\sigma$-finite and positive. We define the class $\mathcal{Z}_\mu$ as $\mathcal{Z}_\mu = \{B \in \mathcal{Z} : \mu(B) < \infty\}$. The symbol $\widehat{N} = \{\widehat{N}(B) : B \in \mathcal{Z}_\mu\}$ indicates a *compensated Poisson random measure* on $(Z, \mathcal{Z})$ with control $\mu$. This means that $\widehat{N}$ is a collection of random variables defined of some probability space $(\Omega, \mathcal{F}, \mathbb{P})$, indexed by the elements of $\mathcal{Z}_\mu$ and such that: (i) for every $B, C \in \mathcal{Z}_\mu$ such that $B \cap C = \varnothing$, $\widehat{N}(B)$ and $\widehat{N}(C)$ are independent; (ii) for every $B \in \mathcal{Z}_\mu$,

$$\widehat{N}(B) \stackrel{\text{law}}{=} \mathfrak{P}(B) - \mu(B),$$

where $\mathfrak{P}(B)$ is a Poisson random variable with parameter $\mu(B)$. Note that properties (i)–(ii) imply, in particular, that $\widehat{N}$ is an *independently scattered* (or *completely random*) measure (see [19]). For every deterministic function $h \in L^2(Z, \mathcal{Z}, \mu) = L^2(\mu)$, we write $\widehat{N}(h) = \int_Z h(z) \widehat{N}(\mathrm{d}z)$ to indicate the Wiener–Itô integral of $h$ with respect to $\widehat{N}$ (see [22]). We recall that for every $h \in L^2(\mu)$, $\widehat{N}(h)$ has an infinitely divisible law, with Lévy–Khinchine exponent (again see [22]) given by

$$\psi(h, \lambda) = \int_Z \exp(\mathrm{i}\lambda h(z) - 1 - \mathrm{i}\lambda h(z)) \mu(\mathrm{d}z), \qquad \lambda \in \mathbb{R}. \tag{7}$$

Also, recall the isometric relation: for every $g, h \in L^2(\mu)$, $\mathbb{E}[\widehat{N}(g)\widehat{N}(h)] = \int_Z h(z)g(z)\mu(\mathrm{d}z)$.

Although our main results involve exclusively single and double integrals with respect to $\widehat{N}$, we will often need to deal with multiple integrals of higher orders. To this end, we briefly recall the definition, as well as some basic properties, of general multiple Poisson integrals. Fix $n \geq 2$. We denote by $L^2(\mu^n)$ the space of real-valued functions on $Z^n$ that are square-integrable with respect to $\mu^n$ and we write $L_s^2(\mu^n)$ to indicate the subspace of $L^2(\mu^n)$ composed of symmetric functions. We denote by $S_n$ and $\widetilde{S}_n$, respectively, the vector space generated by simple functions with the form

$$f(z_1,\ldots,z_n) = \mathbf{1}_{B_1}(z_1)\cdots \mathbf{1}_{B_n}(z_n), \tag{8}$$

where $B_1,\ldots,B_n$ are disjoint subsets of $Z$, and the vector space generated by the symmetrizations of the elements of $S_n$. If $f$ is as in (8) and $\tilde{f} \in \widetilde{S}_n$ is its symmetrization, we define $I_n(\tilde{f})$ as

$$I_n(\tilde{f}) = \widehat{N}(B_1)\cdots \widehat{N}(B_n). \tag{9}$$

The random variable $I_n(\tilde{f})$ is the *multiple Wiener–Itô integral* of order $n$, of $\tilde{f}$ with respect to $\widehat{N}$. Now, note that $\mu$ is non-atomic so $\mu^n$ does not charge diagonals. This implies that for every $n \geq 2$, $\widetilde{S}_n$ is dense in $L_s^2(\mu^n)$ and, consequently, the domain of the operator $I_n$ can be extended to $L_s^2(\mu^n)$ by continuity, due to the isometric formula, true for every $m,n \geq 2$, $\tilde{f} \in \widetilde{S}_n$ and $\tilde{g} \in \widetilde{S}_m$,

$$\mathbb{E}(I_n(\tilde{f})I_m(\tilde{g})) = n!\langle \tilde{f},\tilde{g}\rangle_{L^2(\mu^n)}\mathbf{1}_{(n=m)}. \tag{10}$$

We also use the following conventional notation: $I_1(f) = \widehat{N}(f)$, $f \in L^2(\mu)$; $I_n(f) = I_n(\tilde{f})$, $f \in L^2(\mu^n)$, $n \geq 2$ (this convention extends the definition of $I_n^{\widehat{N}}(f)$ to non-symmetric functions $f$); $I_0(c) = c$, $c \in \mathbb{R}$.

Before stating the main result of the section, we recall a useful version of the *multiplication formula* for multiple Poisson integrals. To this end, we define, for $q,p \geq 1$, $f \in L_s^2(\mu^p)$, $g \in L_s^2(\mu^q)$, $r = 0,\ldots,q \wedge p$ and $l = 1,\ldots,r$, the (contraction) kernel on $Z^{p+q-r-l}$, which reduces the number of variables in the product $fg$ from $p+q$ to $p+q-r-l$ as follows: $r$ variables are identified and, among these, $l$ are integrated out. This contraction kernel is formally defined as follows:

$$f \star_r^l g(\gamma_1,\ldots,\gamma_{r-l},t_1,\ldots,t_{p-r},s_1,\ldots,s_{q-r})$$
$$= \int_{Z^l} f(z_1,\ldots,z_l,\gamma_1,\ldots,\gamma_{r-l},t_1,\ldots,t_{p-r})$$
$$\times g(z_1,\ldots,z_l,\gamma_1,\ldots,\gamma_{r-l},s_1,\ldots,s_{q-r})\mu^l(\mathrm{d}z_1\cdots \mathrm{d}z_l)$$

and, for $l = 0$,

$$\begin{aligned} f \star_r^0 g(\gamma_1,\ldots,\gamma_r,t_1,\ldots,t_{p-r},s_1,\ldots,s_{q-r}) \\ = f(\gamma_1,\ldots,\gamma_r,t_1,\ldots,t_{p-r})g(\gamma_1,\ldots,\gamma_r,s_1,\ldots,s_{q-r}) \end{aligned} \tag{11}$$

so that $f \star_0^0 g(t_1,\ldots,t_p,s_1,\ldots,s_q) = f(t_1,\ldots,t_p)g(s_1,\ldots,s_q)$. For example, if $p = q = 2$,

$$f \star_1^0 g(\gamma,t,s) = f(\gamma,t)g(\gamma,s), \qquad f \star_1^1 g(t,s) = \int_Z f(z,t)g(z,s)\mu(\mathrm{d}z) \tag{12}$$

$$f \star_2^1 g(\gamma) = \int_Z f(z,\gamma)g(z,\gamma)\mu(\mathrm{d}z),$$
$$f \star_2^2 g = \int_Z \int_Z f(z_1,z_2)g(z_1,z_2)\mu(\mathrm{d}z_1)\mu(\mathrm{d}z_2). \tag{13}$$

The following product formula for two Poisson multiple integrals is proved in, for example, [11] and [23]: letting $f \in L_s^2(\mu^p)$ and $g \in L_s^2(\mu^q)$, $p,q \geq 1$, and further supposing that $f \star_r^l g \in L^2(\mu^{p+q-r-l})$ for every $r = 0, \ldots, p \wedge q$ and $l = 1, \ldots, r$, we have further supposing

$$I_p(f)I_q(g) = \sum_{r=0}^{p \wedge q} r! \binom{p}{r}\binom{q}{r} \sum_{l=0}^{r} \binom{r}{l} I_{q+p-r-l}(\widetilde{f \star_r^l g}), \tag{14}$$

where the tilde ($\sim$) stands for symmetrization, that is,

$$\widetilde{f \star_r^l g}(x_1, \ldots, x_{q+p-r-l}) = \frac{1}{(q+p-r-l)!} \sum_\sigma f \star_r^l g(x_{\sigma(1)}, \ldots, x_{\sigma(q+p-r-l)}),$$

where $\sigma$ runs over all $(q+p-r-l)!$ permutations of the set $\{1, \ldots, q+p-r-l\}$.

In what follows, we will systematically work under the assumption that there exists a collection of subsets $\{Z_t : t \in [0,1]\} \subset \mathcal{Z}$ such that $Z_0 = \varnothing$, $Z_1 = Z$, $Z_s \subseteq Z_t$ for $s < t$ and the following continuity condition is satisfied: for every $g \in L^1(\mu)$ and every $t \in [0,1]$,

$$\lim_{s \to t} \int_{Z_s} g(z)\mu(\mathrm{d}z) = \int_{Z_t} g(z)\mu(\mathrm{d}z). \tag{15}$$

Observe that (15) is easily satisfied when $Z$ is a Euclidean space. For instance, when $Z = \mathbb{R}_+^2$ and $\mu$ is equal to the Lebesgue measure, one can take $Z_t = [0, \log(1/1-t)]^2$ ($t \in [0,1)$) and $Z_1 = \mathbb{R}_+^2$. Also, note that for every $t \in [0,1]$, the operator

$$\pi_t : L^2(\mu) \mapsto L^2(\mu) : f \mapsto \pi_t f = f \mathbf{1}_{Z_t} \tag{16}$$

defines a projection operator. Since the family of projections $\pi = \{\pi_t : t \in [0,1]\}$ is non-decreasing and continuous, one says that $\pi$ is a *continuous resolution of the identity* (see, e.g., [28]). This concept in central in the general theory developed in [17]. In what follows, we will use the following notation (introduced in [19]): for every $z, z' \in Z$, one writes

$$z' \prec_\pi z \tag{17}$$

whenever there exists $t \in [0,1]$ such that $z' \in Z_t$ and $z \in Z_t^c$. As an example, consider the case $Z = [0,1]^d$ ($d \geq 1$) and $Z_t = [0,t]^d$, and fix $z = (z^{(1)}, \ldots, z^{(d)}) \in (0,1)^d$. Then $z' \prec_\pi z$ if and only if $z' \in [0, \overline{z}]^d$, where $\overline{z} = \max(z^{(1)}, \ldots, z^{(d)})$. In particular, in the case $d = 1$, one has that $z' \prec_\pi z$ if and only if $z' < z$.

Given a kernel $f \in L_s^2(\mu^2)$ and $z \in Z$, we will write $f\mathbf{1}(\cdot \prec_\pi z)$ to indicate the application

$$Z \mapsto \mathbb{R} : y \mapsto f(y,z)\mathbf{1}(y \prec_\pi z); \tag{18}$$

note that $f\mathbf{1}(\cdot \prec_\pi z) \in L^2(\mu)$. As an example, again consider the case $Z = [0,1]^d$ and $Z_t = [0,t]^d$, $t \in [0,1]$. Then for every fixed $z = (z^{(1)}, \ldots, z^{(d)}) \in (0,1)^d$, the function $f\mathbf{1}(\cdot \prec_\pi z)$ is the application

$$[0,1]^d \mapsto \mathbb{R} : y \mapsto f(y,z)\mathbf{1}(y \in [0,\overline{z}]^d),$$

where $\overline{z} = \max(z^{(1)}, \ldots, z^{(d)})$.

## 3. Main results: CLTs for single and double Poisson integrals

In this section, we apply the results proved in [19] to study the asymptotic behavior of a sequence of random variables of the type

$$F_n = I_2(f_n), \qquad n \geq 1, \tag{19}$$

where $f_n \in L_s^2(\mu^2)$. In particular, our starting point is the following result, taken from [19].

**Proposition 1 (See [19]).** *Consider the sequence $\{F_n : n \geq 1\}$ in (19) and for every fixed $z \in Z$, define $f_n \mathbf{1}(\cdot \prec_\pi z) \in L^2(\mu)$ according to (17) and (18). Suppose that for every $\lambda \in \mathbb{R}$,*

$$\int_Z \exp(\mathrm{i}\lambda 2 \widehat{N}(f_n\mathbf{1}(\cdot \prec_\pi z))) - 1 - \mathrm{i}\lambda 2 \widehat{N}(f_n\mathbf{1}(\cdot \prec_\pi z)))\mu(\mathrm{d}z) \xrightarrow[n\to\infty]{\mathbb{P}} -\frac{\lambda^2}{2}, \tag{20}$$

*where $\xrightarrow{\mathbb{P}}$ denotes convergence in probability. Then $F_n \xrightarrow{\text{law}} \mathcal{N}(0,1)$, where $\mathcal{N}(0,1)$ denotes a standard Gaussian random variable.*

Note that condition (20) is quite difficult to verify, since it involves a continuum of non-trivial transformations of the kernels $f_n$ (one for every $z \in Z$). However, in what follows, we shall show that Proposition 1 can be used to obtain neat sufficient (and sometimes also necessary) conditions on the kernels $f_n$ for the CLT $F_n \xrightarrow{\text{law}} \mathcal{N}(0,1)$ to hold. In particular, these conditions do not involve the resolution of the identity $\pi$ and are exclusively expressed in terms of the kernels $\{f_n\}$. Our results will apply to sequences of kernels satisfying the following assumption.

***Assumption N.*** *The sequence $f_n$, $n \geq 1$, in (19) verifies:*

(N-i) *(integrability condition)* $\forall n \geq 1$,

$$\int_Z f_n(z,\cdot)^2 \mu(\mathrm{d}z) \in L^2(\mu) \quad \text{and} \quad \left\{\int_Z f_n(z,\cdot)^4 \mu(\mathrm{d}z)\right\}^{1/2} \in L^1(\mu); \tag{21}$$

(N-ii) *(normalization condition) as $n \to \infty$,*

$$2 \int_Z \int_Z f_n(z, z')^2 \mu(\mathrm{d}z) \mu(\mathrm{d}z') \to 1; \tag{22}$$

(N-iii) *(fourth power condition) as $n \to \infty$,*

$$\int_Z \int_Z f_n(z, z')^4 \mu(\mathrm{d}z) \mu(\mathrm{d}z') \to 0 \tag{23}$$

*(this implies, in particular, that $f_n \in L^4(\mu^2)$).*

**Remarks.** (1) Suppose there exists a set $B$, independent of $n$, such that $\mu(B) < \infty$ and for each $n$, $f_n = f_n \mathbf{1}_B$, a.e.-$\mathrm{d}\mu^2$ (this is true, in particular, when $\mu$ is finite). Then by the Cauchy–Schwarz inequality, if (23) is verified, $(f_n)$ must necessarily converge to zero in $L_s^2(\mu^2)$. To get more general sequences $(f_n)$, we need to suppose $\mu(Z) = \infty$.

(2) As shown in [13], condition (23) is not required to obtain CLTs for sequences of double Wiener–Itô integrals with respect to Gaussian processes and hence, in this case, it is not necessary to suppose $\mu(Z) = \infty$. Thus, in [13], the setup of a Gaussian measure on $[0,1]^d$ ($d \geq 1$) with a finite control measure was considered.

(3) Assumption N is satisfied, for example, by a properly normalized sequence of uniformly bounded functions, with supports 'slowly converging' to $Z$. For instance, consider a sequence $g_n \in L_s^2(\mu^2)$ such that for $n \geq 1$, $|g_n(\cdot, \cdot\cdot)| \leq c < \infty$ ($c$ independent of $n$) and the support of $g_n$ is contained in a set of the type $B_n \times B_n$, where $0 < \mu(B_n) < \infty$ and $\mu(B_n) \to \infty$. Then if

$$\mu(B_n)^{-2} \int_Z \int_Z g_n(z, z')^2 \mu(\mathrm{d}z) \mu(\mathrm{d}z') \to 1,$$

the sequence $f_n \triangleq \mu(B_n)^{-1} g_n$, $n \geq 1$, verifies Assumption N. Indeed, since $|f_n| \leq c\mu(B_n)^{-1}$,

$$\int_Z \left( \int_Z f_n(z, z')^4 \mu(\mathrm{d}z) \right)^{1/2} \mu(\mathrm{d}z') \leq \frac{c^2}{\mu(B_n)^{1/2}} < \infty,$$

$$\int_Z \int_Z f_n(z, z')^4 \mu(\mathrm{d}z) \mu(\mathrm{d}z') \leq \frac{c^4}{\mu(B_n)^2} \to 0,$$

$$\int_Z \left( \int_Z f_n(z, z')^2 \mu(\mathrm{d}z) \right)^2 \mu(\mathrm{d}z') \leq \frac{c^4}{\mu(B_n)} < \infty.$$

The following central limit theorem is one of the main results of the paper.

**Theorem 2.** *Define the sequence $F_n = I_2(f_n)$ and $f_n \in L_s^2(\mu^2)$, $n \geq 1$, as in (19), and suppose Assumption N holds. Then*

$$f_n \star_1^0 f_n \in L^2(\mu^3) \quad \text{and} \quad f_n \star_1^1 f_n \in L_s^2(\mu^2)$$

*for every $n \geq 1$ and, moreover,*

1. *if*

$$f_n \star_1^1 f_n \to 0 \quad in \ L^2_{s,0}(\mu^2) \quad and \quad f_n \star_2^1 f_n \to 0 \quad in \ L^2(\mu), \quad (24)$$

   *then*

$$F_n \overset{law}{\to} \mathcal{N}(0,1), \quad (25)$$

   *where $\mathcal{N}(0,1)$ is a standard Gaussian random variable;*

2. *if $F_n \in L^4(\mathbb{P})$ for every $n$, then a sufficient condition to have (24) is that*

$$\mathbb{E}(F_n^4) \to 3 = \mathbb{E}[\mathcal{N}(0,1)^4]; \quad (26)$$

3. *if the sequence $\{F_n^4 : n \geq 1\}$ is uniformly integrable, then conditions (24), (25) and (26) are equivalent.*

**Remark.** (a) As already indicated in the Introduction, point 3 can be rephrased by saying that, if Assumption N is satisfied and $\{F_n^4 : n \geq 1\}$ is uniformly integrable, then $F_n \overset{law}{\to} \mathcal{N}(0,1)$ if and only if $\mathbb{E}(F_n^4) \to \mathbb{E}[\mathcal{N}(0,1)^4]$. This result should be compared with the usual 'method of moments' for sequences of random variables. See, for example, [24] and [25].

(b) We recall (see, e.g., [8], page 355) that $\{F_n^4 : n \geq 1\}$ is uniformly integrable if and only if

$$\lim_{M \to \infty} \sup_{n \geq 1} \mathbb{E}[F_n^4 \mathbf{1}_{(F_n^4 > M)}] = 0.$$

(c) While the statement of Theorem 2 does not involve any resolution of the identity $\pi$, the objects appearing in (16) and (17) will play a crucial role in the proof.

(d) Observe that

$$\|f_n \star_1^1 f_n\|^2_{L^2(\mu^2)} = \int_Z \int_Z \left( \int_Z f_n(a,z) f_n(b,z) \mu(\mathrm{d}z) \right)^2 \mu(\mathrm{d}a) \mu(\mathrm{d}b), \quad (27)$$

$$\|f_n \star_2^1 f_n\|^2_{L^2(\mu)} = \int_Z \left( \int_Z f_n(a,z)^2 \mu(\mathrm{d}a) \right)^2 \mu(\mathrm{d}z), \quad (28)$$

$$\|f_n \star_1^0 f_n\|^2_{L^2(\mu^3)} = \int_Z \int_Z \int_Z (f_n(a,z) f_n(b,z))^2 \mu(\mathrm{d}z) \mu(\mathrm{d}a) \mu(\mathrm{d}b), \quad (29)$$

(e) Let $G$ be a Gaussian measure on $(Z, \mathcal{Z})$ with non-atomic control $\mu$ and for $n \geq 1$, let $H_n = I_2^G(h_n)$ be the double Wiener–Itô integral of a function $h_n \in L^2_s(\mu^2)$. In [13] Theorem 1 it is proved that if the normalization relation $2\|h_n\|^2 \to 1$ holds and *regardless of Assumptions (N-i) and (N-iii)*, the following three conditions are equivalent: (i) $H_n \overset{law}{\to} \mathcal{N}(0,1)$; (ii) $\mathbb{E}(H_n^4) \to 3$; (iii) $h_n \star_1^1 h_n \to 0$. Also, note that [13] Theorem 1 applies to multiple integrals of arbitrary order.

(f) A sufficient condition for the uniform integrability of $(F_n^4)$ is clearly that $\sup_n \mathbb{E}(F_n^{4+\varepsilon}) < \infty$ for some $\varepsilon > 0$. Note that in the Gaussian framework of [13] Theorem 1 the uniform integrability condition is always satisfied. Indeed, by noting $H_n = I_2^G(h_n)$ ($n \geq 1$), the sequence of double integrals introduced in the previous remark, for every $p > 2$, there exists a finite constant $c_p$ such that $\sup_n \mathbb{E}(|H_n|^p) \leq c_p (\sup_n \mathbb{E}(H_n^2))^{p/2} < \infty$, where the last relation follows from the normalization condition $\mathbb{E}(H_n^2) = 2\|h_n\|^2 \to 1$.

**Proof of Theorem 2.** Since

$$f_n \star_1^1 f_n(t,s) = \int_Z f_n(s,z) f_n(t,z) \mu(\mathrm{d}z)$$

and $f \in L^2(\mu^2)$, the relation $f_n \star_1^1 f_n \in L_s^2(\mu^2)$ is a consequence of the Cauchy–Schwarz inequality. On the other hand, by (12),

$$\int_{Z^3} (f_n \star_1^0 f_n(\gamma,t,s))^2 \mu^3(\mathrm{d}\gamma,\mathrm{d}t,\mathrm{d}s) = \int_Z \left( \int_Z f_n(\gamma,t)^2 \mu(\mathrm{d}t) \times \int_Z f_n(\gamma,s)^2 \mu(\mathrm{d}s) \right) \mu(\mathrm{d}\gamma)$$

$$= \int_Z \left( \int_Z f_n(\gamma,s)^2 \mu(\mathrm{d}s) \right)^2 \mu(\mathrm{d}\gamma) < \infty,$$

due to part (N-i) Assumption N, so $f_n \star_1^0 f_n \in L^2(\mu^3)$.

(*Proof of point 1.*) According to Proposition 1, to prove point 1, it is sufficient to show that, when Assumption N is verified, condition (24) implies that (20) is satisfied. To this end, we write

$$\int_Z (\exp(\mathrm{i}\lambda 2\widehat{N}(f_n \mathbf{1}(\cdot \prec_\pi z))) - 1 - \mathrm{i}\lambda 2\widehat{N}(f_n \mathbf{1}(\cdot \prec_\pi z))) \mu(\mathrm{d}z)$$

$$= -\frac{\lambda^2}{2} \int_Z (2\widehat{N}(f_n \mathbf{1}(\cdot \prec_\pi z)))^2 \mu(\mathrm{d}z) \qquad (30)$$

$$+ \int_Z \left[ \exp(\mathrm{i}\lambda 2\widehat{N}(f_n \mathbf{1}(\cdot \prec_\pi z))) - 1 \right.$$

$$\left. - \mathrm{i}\lambda 2\widehat{N}(f_n \mathbf{1}(\cdot \prec_\pi z)) + \frac{\lambda^2}{2} (2\widehat{N}(f_n \mathbf{1}(\cdot \prec_\pi z)))^2 \right] \mu(\mathrm{d}z)$$

$$\triangleq U_n + V_n$$

and we shall show that under the assumptions of Theorem 2, $U_n \xrightarrow{\mathbb{P}} -\frac{\lambda^2}{2}$ and $V_n \xrightarrow{\mathbb{P}} 0$. To show that

$$\int_Z (2\widehat{N}(f_n \mathbf{1}(\cdot \prec_\pi z)))^2 \mu(\mathrm{d}z) \xrightarrow{\mathbb{P}} 1 \qquad (31)$$

and hence that $U_n \xrightarrow{\mathbb{P}} -\frac{\lambda^2}{2}$, we apply (14) in the case $p = q = 1$ and obtain

$$((2\widehat{N}(f_n \mathbf{1}(\cdot \prec_\pi z)))^2)^2$$

$$= 4I_1(f_n(z,\cdot)\mathbf{1}(\cdot \prec_\pi z))^2 \tag{32}$$

$$\triangleq 4\int_Z f_n(z,x)^2 \mathbf{1}(x \prec_\pi z)\mu(\mathrm{d}x) + 4I_1(f_n(z,\cdot)^2\mathbf{1}(\cdot \prec_\pi z)) + 4I_2(g_n(z;\cdot,\cdot\cdot)),$$

where $g_n(z;\cdot,\cdot\cdot) \in L^2_s(\mu^2)$ is given by

$$\begin{aligned} g_n(z;a,b) &= f_n \star_0^0 f_n(z,a;z,b)\mathbf{1}_{(a\prec_\pi z)}\mathbf{1}_{(b\prec_\pi z)} \\ &= f_n(z,a)f_n(z,b)\mathbf{1}_{(a\prec_\pi z)}\mathbf{1}_{(b\prec_\pi z)}. \end{aligned} \tag{33}$$

The three terms in (32) correspond, respectively, to the terms $(r=1,l=1)$, $(r=1,l=0)$ and $(r=0,l=0)$ in (14). We deal with each term in (32) in succession. For the first term, observe that, due to [19] Corollary 4 and the symmetry of $f_n$,

$$\int_Z\int_Z f_n(z,x)^2\mu(\mathrm{d}x)\mu(\mathrm{d}z)$$
$$= \int_Z\int_Z f_n(z,x)^2[\mathbf{1}(x \prec_\pi z) + \mathbf{1}(z \prec_\pi x)]\mu(\mathrm{d}x)\mu(\mathrm{d}z)$$
$$= 2\int_Z\int_Z f_n(z,x)^2\mathbf{1}(x \prec_\pi z)\mu(\mathrm{d}x)\mu(\mathrm{d}z)$$

and therefore, thanks to Assumption N,

$$\begin{aligned} -\frac{\lambda^2}{2}\int_Z 4\int_Z f_n(z,x)^2\mathbf{1}(x \prec_\pi z)\mu(\mathrm{d}x)\mu(\mathrm{d}z) \\ = -\frac{\lambda^2}{2}\left[2\int_Z\int_Z f_n(z,x)^2\mu(\mathrm{d}x)\mu(\mathrm{d}z)\right] \to -\frac{\lambda^2}{2}. \end{aligned} \tag{34}$$

For the second term in (32), one has, by applying Fubini's theorem and the Cauchy–Schwarz inequality,

$$\mathbb{E}\left[\int_Z |I_1(f_n(z,\cdot)^2\mathbf{1}(\cdot \prec_\pi z))|\mu(\mathrm{d}z)\right]^2$$
$$= \int_{Z^2} \mathbb{E}[|I_1(f_n(z,\cdot)^2\mathbf{1}(\cdot \prec_\pi z))||I_1(f_n(z',\cdot)^2\mathbf{1}(\cdot \prec_\pi z'))|]\mu(\mathrm{d}z)\mu(\mathrm{d}z')$$
$$\leq \int_{Z^2} \mathbb{E}[I_1(f_n(z,\cdot)^2\mathbf{1}(\cdot \prec_\pi z))^2]^{1/2}\mathbb{E}[I_1(f_n(z',\cdot)^2\mathbf{1}(\cdot \prec_\pi z'))^2]^{1/2}\mu(\mathrm{d}z')\mu(\mathrm{d}z)$$
$$= \left\{\int_Z \left[\int_Z \mu(\mathrm{d}a)f_n(z,a)^4\mathbf{1}(a \prec_\pi z)\right]^{1/2}\mu(\mathrm{d}z)\right\}^2 < \infty$$

since Part ( N–i) of Assumption N is satisfied. It follows that by once again applying Fubini's theorem,

$$\begin{aligned}
&\mathbb{E}\left[\int_Z I_1(f_n(z,\cdot)^2\mathbf{1}(\cdot \prec_\pi z))\mu(\mathrm{d}z)\right]^2 \\
&= \int_{Z^2} \mathbb{E}[I_1(f_n(z,\cdot)^2\mathbf{1}(\cdot \prec_\pi z))I_1(f_n(z',\cdot)^2\mathbf{1}(\cdot \prec_\pi z'))]\mu(\mathrm{d}z)\mu(\mathrm{d}z') \\
&= \int_Z \left[\int_Z f_n(z,x)^2 \mathbf{1}(x \prec_\pi z)\mu(\mathrm{d}z)\right]^2 \mu(\mathrm{d}x) \\
&\leq \int_Z \left[\int_Z f_n(z,x)^2 \mu(\mathrm{d}z)\right]^2 \mu(\mathrm{d}x) \to 0,
\end{aligned} \tag{35}$$

due to (27) and (24). Hence, the integral of the second term in (32) tends to 0 in probability. Now, consider the third term in (32), and observe that by a Fubini argument, and by (33),

$$\begin{aligned}
&\mathbb{E}\left(\int_Z |I_2(g_n(z;\cdot,\cdot\cdot))|\mu(\mathrm{d}z)\right)^2 \\
&= \int_{Z^2} \mathbb{E}(|I_2(g_n(z;\cdot,\cdot\cdot))I_2(g_n(z';\cdot,\cdot\cdot))|)\mu(\mathrm{d}z)\mu(\mathrm{d}z') \\
&\leq \left(\int_Z \mathbb{E}(I_2(g_n(z;\cdot,\cdot\cdot))^2)^{1/2}\mu(\mathrm{d}z)\right)^2 \\
&= \left(\int_Z \left[\int_{Z^2} f_n(z,a)^2 f_n(z,b)^2 \mathbf{1}_{(a\prec_\pi z)} \mathbf{1}_{(b\prec_\pi z)} \mu(\mathrm{d}a)\mu(\mathrm{d}b)\right]^{1/2} \mu(\mathrm{d}z)\right)^2 \\
&= \left(\int_Z \int_Z f_n(z,a)^2 \mathbf{1}_{(a\prec_\pi z)} \mu(\mathrm{d}a)\mu(\mathrm{d}z)\right)^2 < \infty.
\end{aligned}$$

From this, one deduces that

$$\begin{aligned}
&\mathbb{E}\left(\int_Z I_2(g_n(z;\cdot,\cdot\cdot))\mu(\mathrm{d}z)\right)^2 \\
&= \int_{Z^2} \mathbb{E}[I_2(g_n(z;\cdot,\cdot\cdot))I_2(g_n(z';\cdot,\cdot\cdot))]\mu(\mathrm{d}z)\mu(\mathrm{d}z') = 2\|h_n\|^2_{L^2(\mu^2)},
\end{aligned}$$

where, thanks to (33), $h_n \in L^2_s(\mu^2)$ is such that

$$h_n(a,b) = \int_Z f_n(z,a)f_n(z,b)\mathbf{1}_{(a\prec_\pi z)}\mathbf{1}_{(b\prec_\pi z)}\mu(\mathrm{d}z). \tag{36}$$

We now want to show that $f_n \star^1_1 f_n(a,b) = \int_Z f_n(z,a)f_n(z,b)\mu(\mathrm{d}z) \to 0$ implies that $h_n \to 0$. To do this, we start by observing that, a.e.-$\mu^2(\mathrm{d}a,\mathrm{d}b)$ and thanks to Corollary 4 in

[19],
$$f_n(a,b) = f_n(a,b)\mathbf{1}_{(a \prec_\pi b) \cup (b \prec_\pi a)}.$$

As a consequence, by noting (for fixed $z$)

$$(z \prec_\pi a \vee b) = [(z \prec_\pi a) \cap (z \prec_\pi b)] \cup (a \prec_\pi z \prec_\pi b) \cup (b \prec_\pi z \prec_\pi a),$$
$$(a \vee b \prec_\pi z) = (a \prec_\pi z) \cap (b \prec_\pi z),$$
$$(a \vee b \prec_\pi z \vee z') = [(a \vee b \prec_\pi z') \cap (z \prec_\pi z')] \cup [(a \vee b \prec_\pi z) \cap (z' \prec_\pi z)],$$

we obtain that

$$\int_{Z^2} (f_n \star_1^1 f_n(a,b))^2 \mu^2(\mathrm{d}a, \mathrm{d}b)$$
$$= \int_{Z^2} \left( \int_Z f_n(z,a) \mathbf{1}_{(a \prec_\pi z) \cup (z \prec_\pi a)} \mathbf{1}_{(z \prec_\pi b) \cup (b \prec_\pi z)} f_n(z,b) \mu(\mathrm{d}z) \right)^2 \mu^2(\mathrm{d}a, \mathrm{d}b)$$
$$= \int_{Z^2} \left( \int_Z f_n(z,a) f_n(z,b) (\mathbf{1}_{(z \prec_\pi a \vee b)} + \mathbf{1}_{(a \vee b \prec_\pi z)}) \mu(\mathrm{d}z) \right)^2 \mu^2(\mathrm{d}a, \mathrm{d}b) \tag{37}$$
$$= \int_{Z^2} \left( \int_Z f_n(z,a) f_n(z,b) \mathbf{1}_{(z \prec_\pi a \vee b)} \mu(\mathrm{d}z) \right)^2 \mu^2(\mathrm{d}a, \mathrm{d}b) \tag{38}$$
$$+ \int_{Z^2} \left( \int_Z f_n(z,a) f_n(z,b) \mathbf{1}_{(a \vee b \prec_\pi z)} \mu(\mathrm{d}z) \right)^2 \mu^2(\mathrm{d}a, \mathrm{d}b)$$
$$+ 2 \int_{Z^2} \left( \int_Z f_n(z,a) f_n(z,b) \mathbf{1}_{(a \vee b \prec_\pi z)} \mu(\mathrm{d}z) \right)$$
$$\times \left( \int_Z f_n(z',a) f_n(z',b) \mathbf{1}_{(z \prec_\pi a \vee b)} \mu(\mathrm{d}z') \right) \mu^2(\mathrm{d}a, \mathrm{d}b).$$

We now note
$$(a \prec_\pi z \wedge z') = (a \prec_\pi z') \cap (a \prec_\pi z), \tag{39}$$

so that, by a Fubini argument,

$$\int_{Z^2} \left( \int_Z f_n(z,a) f_n(z,b) \mathbf{1}_{(a \vee b \prec_\pi z)} \mu(\mathrm{d}z) \right)^2 \mu^2(\mathrm{d}a, \mathrm{d}b)$$
$$= \int_{Z^2} \left( \int_Z f_n(z,a) f_n(z',a) \mathbf{1}_{(a \prec_\pi z \wedge z')} \mu(\mathrm{d}a) \right)^2 \mu^2(\mathrm{d}z, \mathrm{d}z') \tag{40}$$

and also
$$2 \int_{Z^2} \left( \int_Z f_n(z,a) f_n(z,b) \mathbf{1}_{(a \vee b \prec_\pi z)} \mu(\mathrm{d}z) \right)$$

$$\times \left( \int_Z f_n(z',a)f_n(z',b)\mathbf{1}_{(z'\prec_\pi a\vee b)}\mu(\mathrm{d}z') \right)\mu^2(\mathrm{d}a,\mathrm{d}b)$$

$$= \int_{Z^2} \left( \int_Z f_n(z,a)f_n(z,b)\mathbf{1}_{(a\vee b\prec_\pi z\vee z')}\mu(\mathrm{d}z) \right)$$

$$\times \left( \int_Z f_n(z',a)f_n(z',b)\mathbf{1}_{(z'\wedge z\prec_\pi a\vee b)}\mu(\mathrm{d}z') \right)\mu^2(\mathrm{d}a,\mathrm{d}b),$$

so that the relation

$$\mathbf{1}_{(a\vee b\prec_\pi z\vee z')}\mathbf{1}_{(z'\wedge z\prec_\pi a\vee b)} = \mathbf{1}_{(z'\wedge z\prec_\pi a, b\prec_\pi z\vee z')} + \mathbf{1}_{(z'\wedge z\prec_\pi a\prec_\pi z\vee z')}\mathbf{1}_{(b\prec_\pi z\wedge z')}$$

$$+ \mathbf{1}_{(z'\wedge z\prec_\pi b\prec_\pi z\vee z')}\mathbf{1}_{(a\prec_\pi z\wedge z')}$$

gives

$$2\int_{Z^2} \left( \int_Z f_n(z,a)f_n(z,b)\mathbf{1}_{(a\vee b\prec_\pi z)}\mu(\mathrm{d}z) \right)$$

$$\times \left( \int_Z f_n(z',a)f_n(z',b)\mathbf{1}_{(z'\prec_\pi a\vee b)}\mu(\mathrm{d}z') \right)\mu^2(\mathrm{d}a,\mathrm{d}b)$$

$$= \int_{Z^2} \left( \int_Z f_n(z,a)f_n(z',a)\mathbf{1}_{(z\wedge z'\prec_\pi a\prec_\pi z\vee z')}\mu(\mathrm{d}a) \right)^2 \mu^2(\mathrm{d}z,\mathrm{d}z') \tag{41}$$

$$+ 2\int_{Z^2} \left( \int_Z f_n(z,a)f_n(z',a)\mathbf{1}_{(z\wedge z'\prec_\pi a\prec_\pi z\vee z')}\mu(\mathrm{d}a) \right) \tag{42}$$

$$\times \left( \int_Z f_n(z,a)f_n(z',a)\mathbf{1}_{(a\prec_\pi z\wedge z')}\mu(\mathrm{d}a) \right)\mu^2(\mathrm{d}z,\mathrm{d}z').$$

After making the change of variables $(z,a,b) \to (a,z,z')$ in (38), observe that the terms (38), (40), (41) and (42) are integrals of terms of the form $(A+B)^2$, $A^2$, $B^2$ and $2AB$, respectively, whose sum therefore equals $2(A+B)^2$, yielding

$$\int_{Z^2} (f_n \star_1^1 f_n(a,b))^2 \mu^2(\mathrm{d}a,\mathrm{d}b)$$

$$= 2\int_{Z^2} \left( \int_Z f_n(z,a)f_n(z,b)\mathbf{1}_{(z\prec_\pi a\vee b)}\mu(\mathrm{d}z) \right)^2 \mu^2(\mathrm{d}a,\mathrm{d}b). \tag{43}$$

Since $h_n$ (as defined in (36)) is such that

$$\int_{Z^2} h_n(a,b)^2 \mu^2(\mathrm{d}a,\mathrm{d}b) = \int_{Z^2} \left( \int_Z f_n(z,a)f_n(z,b)\mathbf{1}_{(a\vee b\prec_\pi z)}\mu(\mathrm{d}z) \right)^2 \mu^2(\mathrm{d}a,\mathrm{d}b)$$

and

$$\int_{Z^2} (f_n \star_1^1 f_n(a,b))^2 \mu^2(\mathrm{d}a,\mathrm{d}b)$$

$$= \int_{Z^2} \left[ \int_Z f_n(z,a) f_n(z,b) (\mathbf{1}_{(a \vee b \prec_\pi z)} + \mathbf{1}_{(z \prec_\pi a \vee b)}) \mu(dz) \right]^2 \mu^2(da, db),$$

relation (43) gives the implication: if $f_n \star_1^1 f_n \to 0$ in $L^2(\mu^2)$, then

$$h_n \to 0 \quad \text{in } L^2(\mu^2). \tag{44}$$

This last result, combined with (34) and (35), implies that the sequence $U_n$, $n \geq 1$, as defined in (30), converges to $-\frac{\lambda^2}{2}$ in probability.

To show that $V_n \xrightarrow{\mathbb{P}} 0$, observe that $|\exp(i\lambda x) - 1 - i\lambda x + \frac{1}{2}\lambda^2 x^2| \leq |\lambda x|^3/6$ and, consequently, by Cauchy–Schwarz,

$$|V_n| \leq \frac{|\lambda|^3}{6} \int_Z |2\widehat{N}(f_n \mathbf{1}(\cdot \prec_\pi z))|^3 \mu(dz)$$

$$\leq \frac{|\lambda|^3}{6} \left( \int_Z |2\widehat{N}(f_n \mathbf{1}(\cdot \prec_\pi z))|^4 \mu(dz) \right)^{1/2} \left( \int_Z |2\widehat{N}(f_n \mathbf{1}(\cdot \prec_\pi z))|^2 \mu(dz) \right)^{1/2}. \tag{45}$$

Since the first part of the proof implies that under (24), $(\int_Z |2\widehat{N}(f_n \mathbf{1}(\cdot \prec_\pi z))|^2 \mu(dz))^{1/2} \xrightarrow{\mathbb{P}} 1$, to conclude the proof of point 1 it is sufficient to show that under Assumption N and (24),

$$\int_Z |2\widehat{N}(f_n \mathbf{1}(\cdot \prec_\pi z))|^4 \mu(dz) \to 0$$

in $L^1(\mathbb{P})$. To do this, one can use (32) and the orthogonality of multiple integrals of different orders to obtain that for any fixed $z$,

$$\mathbb{E}[(2\widehat{N}(f_n \mathbf{1}(\cdot \prec_\pi z)))^4]$$
$$= \mathbb{E}[((2\widehat{N}(f_n \mathbf{1}(\cdot \prec_\pi z)))^2)^2]$$
$$= 16 \left( \int_Z f_n(z,x)^2 \mathbf{1}(x \prec_\pi z) \mu(dx) \right)^2 + 16 \int_Z f_n(z,\cdot)^4 \mathbf{1}(\cdot \prec_\pi z) \mu(da)$$
$$+ 32 \int_{Z^2} f_n(z,a)^2 f_n(z,b)^2 \mathbf{1}_{(a \prec_\pi z)} \mathbf{1}_{(b \prec_\pi z)} \mu^2(da, db)$$

and therefore

$$\mathbb{E} \int_Z |2\widehat{N}(f_n \mathbf{1}(\cdot \prec_\pi z))|^4 \mu(dz)$$
$$= \int_Z \mathbb{E} |2\widehat{N}(f_n \mathbf{1}(\cdot \prec_\pi z))|^4 \mu(dz)$$
$$\leq 16 \int_Z \left( \int_Z f_n(z,x)^2 \mu(dx) \right)^2 \mu(dz)$$

$$+ 16 \int_Z \int_Z f_n(z,a)^4 \mu(\mathrm{d}a)\mu(\mathrm{d}z) + 32 \int_Z \left( \int_Z f_n(z,x)^2 \mu(\mathrm{d}x) \right)^2 \mu(\mathrm{d}z)$$
$$\to 0,$$

since Assumption N and (24) are in order. This concludes the proof of part 1.

(*Proof of point 2.*) To proof point 2, we use the product formula expansion (14) (from the term with $r=0$ to the terms with $r=2$) to write

$$F_n^2 = I_3(f_n)^2 = I_4(\widetilde{f_n \star_0^0 f_n}) + 4I_3(\widetilde{f_n \star_1^0 f_n}) + 4I_2(f_n \star_1^1 f_n)$$
$$+ 2I_2(f_n \star_2^0 f_n) + 2I_1(f_n \star_2^1 f_n) + 2\|f_n\|_{L_s^2(\mu^2)}^2$$

and observe that since Assumption N holds and $f_n \star_2^0 f_n(a,b) = f_n(a,b)^2$ (by (11)), $I_2(f_n \star_2^0 f_n) \to 0$ in $L^2(\mathbb{P})$ by (23) and therefore the assumption $\mathbb{E}(F_n^4) \to 3$ implies that

$$\mathbb{E}[(F_n^2 - 2I_2(f_n \star_2^0 f_n))^2] \tag{46}$$
$$= \mathbb{E}[(I_4(\widetilde{f_n \star_0^0 f_n}) + 4I_3(\widetilde{f_n \star_1^0 f_n}) + 4I_2(f_n \star_1^1 f_n) + 2I_1(f_n \star_2^1 f_n) + 2\|f_n\|_{L_s^2(\mu^2)}^2)^2]$$
$$\to 3. \tag{47}$$

Now, due to (46),

$$\mathbb{E}[(F_n^2 - 2I_2(f_n \star_2^0 f_n))^2] \tag{48}$$
$$= \mathbb{E}(I_4(\widetilde{f_n \star_0^0 f_n})^2) + \mathbb{E}(16 I_3(\widetilde{f_n \star_1^0 f_n})^2) + \mathbb{E}(16 I_2(f_n \star_1^1 f_n)^2) \tag{49}$$
$$+ \mathbb{E}(4 I_1(f_n \star_2^1 f_n)^2) + (2\|f_n\|_{L_s^2(\mu^2)}^2)^2.$$

There are no cross terms because the multiple integrals have different orders and hence are orthogonal. The most complicated term in the square bracket in (49) is $\mathbb{E}(I_4(\widetilde{f_n \star_0^0 f_n})^2)$. Since we are dealing with second order moments, the computations are as in the Gaussian case. We therefore obtain, using, for example, [20], formula (2), page 250, that

$$\mathbb{E}(I_4(\widetilde{f_n \star_0^0 f_n})^2) = 4! \|\widetilde{f_n \star_0^0 f_n}\|_{L^2(\mu^2)}^2 = 2(2\|f_n\|_{L_s^2(\mu^2)}^2)^2 + 16\|f_n \star_1^1 f_n\|_{L^2(\mu^2)}^2.$$

As a consequence, (48) equals

$$[2(2\|f_n\|_{L_{s,0}^2(\mu^2)}^2)^2 + 16\|f_n \star_1^1 f_n\|_{L^2(\mu^2)}^2] + 16 \times 3! \|f_n \star_1^0 f_n\|_{L^2(\mu^3)}^2$$
$$+ 16 \times 2\|f_n \star_1^1 f_n\|_{L^2(\mu^2)}^2 + 4\|f_n \star_2^1 f_n\|_{L^2(\mu)}^2 + (2\|f_n\|_{L_s^2(\mu^2)}^2)^2 \tag{50}$$
$$= 3(2\|f_n\|_{L_s^2(\mu^2)}^2)^2 + 48\|f_n \star_1^1 f_n\|_{L^2(\mu^2)}^2 + 96\|f_n \star_1^0 f_n\|_{L^2(\mu^3)}^2 + 4\|f_n \star_2^1 f_n\|_{L^2(\mu)}^2.$$

Since (50) converges to 3 by (47) and $2\|f_n\|_{L_s^2(\mu^2)}^2 \to 1$ by Assumption (N-ii), we conclude that $\|f_n \star_1^1 f_n\|_{L^2(\mu^2)}^2 \to 0$ and $\|f_n \star_2^1 f_n\|_{L^2(\mu)}^2 \to 0$, thus proving point 2.

(*Proof of point 3.*) If $F_n \overset{\text{law}}{\to} \mathcal{N}(0,1)$ and $(F_n^4)$ is uniformly integrable, then necessarily $\mathbb{E}(F_n^4) \to \mathbb{E}(\mathcal{N}(0,1)^4) = 3$. We had $(26) \Rightarrow (25) \Rightarrow (24)$ and we just showed that uniform integrability and $(24)$ imply $(26)$, proving the equivalence of the three statements under uniform integrability. $\square$

***Example.*** We exhibit an elementary example of a sequence $f_n \in L_s^2(\mu^2)$, $n \geq 1$, verifying conditions $(21)$, $(22)$, $(23)$ and $(24)$. As discussed below, this example involves a sequence of i.i.d. random variables of the type $\widehat{N}(B)$. As such, it could be alternatively worked out by means of a standard application of the usual CLT for i.i.d. random variables with second moments. However, it provides a first and simple illustration of our techniques, also showing that our results are consistent with the classic limit results of probability theory. More sophisticated examples are discussed in Section 4. Let $B_j$, $j \geq 1$, be a sequence of disjoint subsets of $Z$ such that $\mu(B_j) = 1$, $j \geq 1$, and set

$$B_{0,j}^2 = \{(x,y) \in B_j \times B_j : x \neq y\}, \qquad j \geq 1.$$

Note that since $\mu$ is non-atomic, $\mu^2(B_{0,j}^2) = \mu^2(B_j \times B_j) = 1$. For $n \geq 1$ and $(x,y) \in Z^2$, we define $f_n(x,y) = (2n)^{-1/2} \sum_{j=1}^n \mathbf{1}_{B_{0,j}^2}(x,y)$. Of course, $f_n \in L_s^2(\mu^2)$, by definition, and we shall prove that $(f_n)$ also satisfies $(21)$, $(22)$, $(23)$ and $(24)$. Indeed, $\int_Z f_n(z,\cdot)^2 \mu(\mathrm{d}z) = 2n^{-1} \sum_{j=1}^n \mathbf{1}_{B_j}(\cdot) \in L^2(\mu)$ and

$$2\|f_n\|^2 = 2(2n)^{-1} \sum_{j=1}^n \mu^2(B_{0,j}^2) = 1,$$

so $(f_n)$ satisfies $(21)$ and $(22)$. On the other hand,

$$\int_Z \int_Z f_n(x,y)^4 \mu^2(\mathrm{d}x,\mathrm{d}y) = \frac{1}{4n^2} \int_Z \int_Z \left(\sum_{j=1}^n \mathbf{1}_{B_{0,j}^2}(x,y)\right)^4 \mu^2(\mathrm{d}x,\mathrm{d}y) = \frac{1}{4n} \to 0$$

and therefore $(23)$ is satisfied. Finally,

$$\int_Z \left(\int_Z f_n(z,x)^4 \mu(\mathrm{d}z)\right)^{1/2} \mu(\mathrm{d}x) \leq \sqrt{n}/2 < \infty, \tag{51}$$

$$\int_Z \left(\int_Z f_n(z,x)^2 \mu(\mathrm{d}z)\right)^2 \mu(\mathrm{d}x) = 1/4n \to 0 \tag{52}$$

and

$$\int_Z \left(\int_Z f_n(x,z) f_n(y,z) \mu(\mathrm{d}z)\right)^2 \mu^2(\mathrm{d}x,\mathrm{d}y)$$
$$= \frac{1}{4n^2} \int_Z \int_Z \left(\sum_{j=1}^n \mathbf{1}_{B_{0,j}^2}(x,y)\right)^2 \mu^2(\mathrm{d}x,\mathrm{d}y) = \frac{1}{4n} \to 0,$$

thus yielding that $(f_n)$ satisfies (24), by (27) and (28). To see (51), recall that $\mu(B_j) = 1$ and write

$$\int_Z \left(\int_Z f_n(z,x)^4 \mu(\mathrm{d}z)\right)^{1/2} \mu(\mathrm{d}x) \leq \frac{1}{2n}\mu\left(\bigcup_{j=1}^n B_j\right)^{1/2} \int_Z \left(\sum_{j=1}^n \mathbf{1}_{B_j}(x)\right)^{1/2} \mu(\mathrm{d}x)$$

$$= \frac{1}{2n}\mu\left(\bigcup_{j=1}^n B_j\right)^{3/2}$$

(we used the fact that $(\sum_{j=1}^n \mathbf{1}_{B_j}(x))^{1/2} = \sum_{j=1}^n \mathbf{1}_{B_j}(x)$). As anticipated, since (due to, e.g., (14))

$$I_2(f_n) = n^{-1/2} \sum_{j=1}^n 2^{-1/2}(\widehat{N}(B_j)^2 - \widehat{N}(B_j) - 1),$$

the central limit result $I_2(f_n) \overset{\text{law}}{\to} \mathcal{N}(0,1)$ can be verified directly by using a standard version of the central limit theorem, as well as the fact that $\widehat{N}$ is independently scattered and the $B_j$'s are disjoint, with $\mu(B_j) = 1$.

By combining the results of [19] and Theorem 2, we can also characterize the joint (weak) convergence of a single and a double Poisson integral toward a bivariate Gaussian vector.

**Theorem 3.** (A) *Consider a sequence*

$$G_n = I_1(g_n) = \int_Z g_n(z)\widehat{N}(\mathrm{d}z), \qquad n \geq 1, \tag{53}$$

*where $g_n \in L^2(\mu) \cap L^3(\mu)$, and suppose that as $n \to \infty$,*

$$\|g_n\|_{L^2(\mu)}^2 \to 1 \quad \text{and} \quad \int_Z |g_n(z)|^3 \mu(\mathrm{d}z) \to 0. \tag{54}$$

*Then $G_n \overset{\text{law}}{\to} X$, where $X$ is a centered standard Gaussian random variable.*

(B) *Consider a sequence $F_n = I_2(f_n)$, $n \geq 1$, with $f_n \in L^2_{s,0}(\mu^2)$ as in (19), and a sequence $G_n = I_1(g_n)$, $n \geq 1$, as in part* (A) *above. Suppose, moreover, that:*

(i) *the sequence $(f_n)$ verifies Assumption N and satisfies condition (24);*
(ii) *the sequence $(g_n)$ satisfies (54).*

*Then as $n \to \infty$,*

$$(F_n, G_n) \overset{\text{law}}{\to} (X, X'), \tag{55}$$

*where $X, X'$ are two independent, centered standard Gaussian random variables.*

***Remark.*** Roughly speaking, Theorem 3 tells us that for sequences of random variables such as $(F_n)$ and $(G_n)$, the (componentwise) weak convergence of $F_n$ and $G_n$ toward a Gaussian distribution *always implies* the joint convergence of the vector $(F_n, G_n)$. As proved in [20], an analogous property holds for vectors of multiple Wiener–Itô integrals with respect to general Gaussian processes.

**Proof of Theorem 3.** (A) By using (7), for every $\lambda \in \mathbb{R}$,

$$\mathbb{E}[\exp(i\lambda G_n)]$$
$$= \exp\left\{\int_Z [\exp(i\lambda g_n(z)) - 1 - i\lambda g_n(z)]\mu(dz)\right\}$$
$$= \exp\left(-\frac{\lambda^2}{2}\|g_n\|^2_{L^2(\mu)}\right) \exp\left(\int_Z \left[\exp(i\lambda g_n(z)) - 1 - i\lambda g_n(z) + \frac{\lambda^2}{2}g_n(z)^2\right]\mu(dz)\right).$$

To conclude, observe that by (54), $\|g_n\|^2_{L^2(\mu)} \to 1$ and the integral in the second exponential is bounded by

$$\int_Z \left|\exp(i\lambda g_n(z)) - 1 - i\lambda g_n(z) + \frac{\lambda^2}{2}g_n(z)^2\right|\mu(dz) \leq \frac{|\lambda|^3}{6}\int_Z |g_n(z)|^3 \mu(dz) \to 0,$$

thus yielding $\mathbb{E}[\exp(i\lambda G_n)] \to \exp(-\lambda^2/2)$, as required.

(B) To prove (55), it is sufficient to show that for every $\alpha, \beta \in \mathbb{R}$,

$$T_n^{(\alpha,\beta)} \triangleq \alpha F_n + \beta G_n \overset{\text{law}}{\to} \alpha X + \beta X' \overset{\text{law}}{=} \sqrt{\alpha^2 + \beta^2} \times X. \tag{56}$$

We shall prove relation (56) by using [19] Theorem 7, as well as some estimates appearing in the proof of Theorem 2. To do this, define $X_{\widehat{N}}(h) \triangleq \int_Z h(z)\widehat{N}(dz) = I_1(h)$ for every $h \in L^2(Z, \mathcal{Z}, \mu) \triangleq \mathfrak{H}_\mu$ and consider the random field

$$X_{\widehat{N}} = \{X_{\widehat{N}}(h) : h \in \mathfrak{H}_\mu\}.$$

Note that $X_{\widehat{N}}$ belongs to the class of random fields studied in [19]. Now, for every $n$, denote by $\pi$ the resolution of the identity over the Hilbert space $\mathfrak{H}_\mu$ given by

$$\pi_t h(z) = \mathbf{1}_{Z_t}(z)h(z), \qquad h \in \mathfrak{H}_\mu, t \in [0,1],$$

where the sets $Z_t$ appear in (15), and further define

$$u_n^{(\alpha,\beta)}(z) = 2\alpha I_1(f_n(\cdot, z)\mathbf{1}(\cdot \prec_\pi z)) + \beta g_n(z)$$
$$= \alpha h_\pi(f_n)(z) + \beta g_n(z) \in L^2_\pi(\mathfrak{H}_\mu, X_{\widehat{N}}),$$

where the class $L^2_\pi(\mathfrak{H}_\mu, X_{\widehat{N}})$ of $\pi$-adapted random functions is defined in Section 3.2 of [17] and the process $z \mapsto h_\pi(f_n)(z) \in L^2_\pi(\mathfrak{H}_\mu, X_{\widehat{N}})$ is defined as $h_\pi(f_n)(z) = 2I_1(f_n(\cdot, z)\mathbf{1}(\cdot \prec_\pi$

$z))$. By using [19] Proposition 5, we can write

$$T_n^{(\alpha,\beta)} = J_{X_{\widehat{N}}}^{\pi}(u_n^{(\alpha,\beta)}) = \alpha J_{X_{\widehat{N}}}^{\pi}(h_\pi(f_n)) + \beta J_{X_{\widehat{N}}}^{\pi}(g_n),$$

where the generalized stochastic integral $J_{X_{\widehat{N}}}^{\pi}$ is defined in Section 3.2 of [17] and the second equality is a consequence of the linearity of the operator $J_{X_{\widehat{N}}}^{\pi}$. From [17] Theorem 7, in we deduce that (56) is proved, once it is shown that

$$\int_Z [\exp(\mathrm{i}\lambda u_n^{(\alpha,\beta)}(z)) - 1 - \mathrm{i}\lambda u_n^{(\alpha,\beta)}(z)]\mu(\mathrm{d}z) \xrightarrow{\mathbb{P}} -\frac{\lambda^2}{2}(\alpha^2 + \beta^2). \tag{57}$$

To prove (57), we adopt the same strategy as in the proof of Theorem 2, that is, we write

$$\int_Z [\exp(\mathrm{i}\lambda u_n^{(\alpha,\beta)}(z)) - 1 - \mathrm{i}\lambda u_n^{(\alpha,\beta)}(z)]\mu(\mathrm{d}z)$$

$$= -\frac{\lambda^2}{2}\int_Z (u_n^{(\alpha,\beta)}(z))^2 \mu(\mathrm{d}z)$$

$$+ \int_Z \left[\exp(\mathrm{i}\lambda u_n^{(\alpha,\beta)}(z)) - 1 - \mathrm{i}\lambda u_n^{(\alpha,\beta)}(z) + \frac{\lambda^2}{2}(u_n^{(\alpha,\beta)}(z))^2\right]\mu(\mathrm{d}z)$$

$$\triangleq A_n + B_n$$

and show that under Assumptions (i)–(iii) in the statement, $A_n \xrightarrow{\mathbb{P}} -\frac{\lambda^2}{2}(\alpha^2 + \beta^2)$ and $B_n \xrightarrow{\mathbb{P}} 0$. We start by considering $B_n$, writing

$$|B_n|^{1/3} = \left|\int_Z \left[\exp(\mathrm{i}\lambda u_n^{(\alpha,\beta)}(z)) - 1 - \mathrm{i}\lambda u_n^{(\alpha,\beta)}(z) + \frac{\lambda^2}{2}(u_n^{(\alpha,\beta)}(z))^2\right]\mu(\mathrm{d}z)\right|^{1/3}$$

$$\leq \frac{|\lambda|}{6^{1/3}}\left(\int_Z |u_n^{(\alpha,\beta)}(z)|^3 \mu(\mathrm{d}z)\right)^{1/3} = \frac{|\lambda|}{6^{1/3}}\left(\int_Z |\alpha h_\pi(f_n)(z) + \beta g_n(z)|^3 \mu(\mathrm{d}z)\right)^{1/3}$$

$$\leq \frac{|\alpha\lambda|}{6^{1/3}}\left(\int_Z |h_\pi(f_n)(z)|^3 \mu(\mathrm{d}z)\right)^{1/3} + \frac{|\beta\lambda|}{6^{1/3}}\left(\int_Z |g_n(z)|^3 \mu(\mathrm{d}z)\right)^{1/3}.$$

Now recall that in the proof of Theorem 2, we already verified in (45) that under Assumption N and (24),

$$\int_Z |h_\pi(f_n)(z)|^3 \mu(\mathrm{d}z) \xrightarrow{\mathbb{P}} 0.$$

Since $(\int_Z |g_n(z)|^3 \mu(\mathrm{d}z))^{1/3} \to 0$ (due to (54)), we deduce that $B_n \xrightarrow{\mathbb{P}} 0$. To prove the convergence of the sequence $A_n$, we write

$$\int_Z (u_n^{(\alpha,\beta)}(z))^2 \mu(\mathrm{d}z)$$

$$= \alpha^2 \int_Z h_\pi(f_n)(z)^2 \mu(\mathrm{d}z) + \beta^2 \int_Z g_n(z)^2 \mu(\mathrm{d}z)$$

$$+ 2\alpha\beta \int_Z h_\pi(f_n)(z)g_n(z)\mu(\mathrm{d}z).$$

In the proof of Theorem 2, we showed that when Assumption N and (24) hold, $\int_Z h_\pi(f_n)(z)^2 \times \mu(\mathrm{d}z) \xrightarrow{\mathbb{P}} 1$ (see (31)). Since $\int_Z g_n(z)^2 \mu(\mathrm{d}z) \to 1$ by (54), to prove that $A_n \xrightarrow{\mathbb{P}} -\frac{\lambda^2}{2}(\alpha^2 + \beta^2)$, we have only to show that

$$\int_Z h_\pi(f_n)(z) g_n(z) \mu(\mathrm{d}z) \xrightarrow{\mathbb{P}} 0.$$

To this end, we use the definition of $h_\pi(f_n)$ and apply Fubini's theorem to obtain that

$$\mathbb{E}\left[\int_Z |h_\pi(f_n)(z)g_n(z)|\mu(\mathrm{d}z)\right]^2$$

$$\leq \int_{Z^2} \left\{\int_Z f_n(a,z)^2 \mathbf{1}(a \prec_\pi z) \mu(\mathrm{d}a) \int_Z f_n(b,z')^2 \mathbf{1}(b \prec_\pi z') \mu(\mathrm{d}b)\right\}^{1/2}$$

$$\times |g_n(z')g_n(z)|\mu(\mathrm{d}z)\mu(\mathrm{d}z')$$

$$= \left(\int_Z |g_n(z)| \left\{\int_Z f_n(b,z)^2 \mathbf{1}(b \prec_\pi z) \mu(\mathrm{d}b)\right\}^{1/2} \mu(\mathrm{d}z)\right)^2$$

$$\leq \int_Z g_n(z)^2 \mu(\mathrm{d}z) \times \int_Z \left\{\int_Z f_n(b,z')^2 \mathbf{1}(b \prec_\pi z') \mu(\mathrm{d}b)\right\} \mu(\mathrm{d}z') < \infty.$$

As a consequence, by once again using the isometric properties of the random measure $\widehat{N}$, Fubini's theorem and Cauchy–Schwarz,

$$\tfrac{1}{4}\mathbb{E}\left[\left(\int_Z h_\pi(f_n)(z)g_n(z)\mu(\mathrm{d}z)\right)^2\right]$$

$$= \int_{Z^2} \mathbb{E}[I_1(f_n(\cdot,z)\mathbf{1}(\cdot \prec_\pi z))g_n(z)I_1(f_n(\cdot,z')\mathbf{1}(\cdot \prec_\pi z'))g_n(z')]\mu(\mathrm{d}z)\mu(\mathrm{d}z')$$

$$= \int_Z \left(\int_Z f_n(a,z)\mathbf{1}(a \prec_\pi z)g_n(z)\mu(\mathrm{d}z)\right)^2 \mu(\mathrm{d}a)$$

$$= \int_Z \left(\int_Z f_n(a,z)\mathbf{1}(a \prec_\pi z)g_n(z)\mu(\mathrm{d}z)\right)$$

$$\times \left(\int_Z f_n(a,z')\mathbf{1}(a \prec_\pi z')g_n(z')\mu(\mathrm{d}z')\right) \mu(\mathrm{d}a)$$

$$= \int_{Z^2} g_n(z)g_n(z') \left(\int_Z f_n(a,z)f_n(a,z')\mathbf{1}(a \prec_\pi z \wedge z')\right) \mu^2(\mathrm{d}z,\mathrm{d}z')$$

$$\leq \left( \int_{Z^2} g_n(z)^2 g_n(z')^2 \mu^2(\mathrm{d}z, \mathrm{d}z') \right)^{1/2}$$
$$\times \left( \int_{Z^2} \left( \int_Z f_n(a,z) f_n(a,z') \mathbf{1}(a \prec_\pi z \wedge z') \mu(\mathrm{d}a) \right)^2 \mu^2(\mathrm{d}z, \mathrm{d}z') \right)^{1/2}.$$

Note that $(\int_{Z^2} g_n(z)^2 g_n(z')^2 \mu^2(\mathrm{d}z, \mathrm{d}z'))^{1/2} = \int_Z g_n(z)^2 \mu(\mathrm{d}z) \to 1$, by (54). Moreover,

$$\int_{Z^2} \left( \int_Z f_n(a,z) f_n(a,z') \mathbf{1}(a \prec_\pi z \wedge z') \mu(\mathrm{d}a) \right)^2 \mu^2(\mathrm{d}z, \mathrm{d}z')$$
$$= \int_{Z^2} \left( \int_Z f_n(a,z) f_n(b,z) \mathbf{1}(a \vee b \prec_\pi z) \mu(\mathrm{d}z) \right)^2 \mu^2(\mathrm{d}a, \mathrm{d}b)$$
$$= \int_{Z^2} h_n(a,b)^2 \mu^2(\mathrm{d}a, \mathrm{d}b) \to 0,$$

where we again applied Fubini's theorem, $h_n$ is given by (36) and we used (44). It follows that $A_n \xrightarrow{\mathbb{P}} -\frac{\lambda^2}{2}(\alpha^2 + \beta^2)$ and therefore that (56) holds, thus concluding the proof of Theorem 3. $\square$

## 4. Applications

This section contains two illustrations of our techniques. Both involve *generalized Volterra processes*, that is, random processes having the form (5), where $\widehat{N}$ is a Poisson measure over $\mathbb{R} \times \mathbb{R}$ (e.g., with control measure $\nu(\mathrm{d}u)\,\mathrm{d}s$) and $h$ is a deterministic bivariate kernel. In Section 4.1, we prove CLTs involving linear and quadratic functionals of Ornstein–Uhlenbeck Lévy processes, obtained from (5) by setting

$$h(t,s) = \sqrt{2\lambda} \exp(-\lambda(t-s)) \mathbf{1}_{t>s}, \qquad \lambda > 0.$$

See [3], [10] and [26], respectively, for applications of Ornstein–Uhlenbeck Lévy processes to finance, survival analysis and network modeling. Section 4.2 contains a concise description of the applications of the theory of this paper developed in [5] and [15], where Theorem 2 and Theorem 3 are applied in order to obtain CLTs involving random hazard rates in Bayesian survival analysis.

### 4.1. Ornstein–Uhlenbeck Lévy processes

Fix $\lambda > 0$. We consider the *Ornstein–Uhlenbeck Lévy process* given by

$$Y_t^\lambda = \sqrt{2\lambda} \int_{-\infty}^t \int_{\mathbb{R}} u \exp(-\lambda(t-x)) \widehat{N}(\mathrm{d}u, \mathrm{d}x), \qquad t \geq 0, \tag{58}$$

where $\widehat{N}$ is a centered Poisson measure over $\mathbb{R} \times \mathbb{R}$ with control measure given by $\nu(\mathrm{d}u)\,\mathrm{d}x$, where $\nu(\cdot)$ is positive and normalized in such a way that $\int_{\mathbb{R}} u^2 \nu(\mathrm{d}u) = 1$. We also assume that $\int_{\mathbb{R}} |u|^3 \nu(\mathrm{d}u) < \infty$. Note that $Y_t^\lambda$ is a stationary moving average Lévy process. We shall use part A of Theorem 3 to prove the following result involving linear functionals of $Y^\lambda$.

**Theorem 4.** *As $T \to \infty$,*

$$\frac{1}{\sqrt{T}} \int_0^T Y_t^\lambda \,\mathrm{d}t \stackrel{\text{law}}{\to} \mathcal{N}(0, \sigma^2(\lambda)), \tag{59}$$

*where $\mathcal{N}(0, \sigma^2(\lambda))$ stands for a centered Gaussian random variable with variance $\sigma^2(\lambda) = 2\lambda^{-1}$.*

**Proof.** First, put the integral in the form (53). Applying Fubini's theorem, we have

$$\frac{1}{\sqrt{T}} \int_0^T Y_t^\lambda \,\mathrm{d}t$$
$$= \int_{-\infty}^T \int_{\mathbb{R}} u \left[ \left( \frac{2\lambda}{T} \right)^{1/2} \int_{x \vee 0}^T \exp(-\lambda(t-x)) \,\mathrm{d}t \right] \widehat{N}(\mathrm{d}u, \mathrm{d}x).$$

We need to verify (54). First

$$\mathbb{E}\left[ \left( \frac{1}{\sqrt{T}} \int_0^T Y_t^\lambda \,\mathrm{d}t \right)^2 \right]$$
$$= \int_{-\infty}^0 \int_{\mathbb{R}} u^2 \left[ \left( \frac{2\lambda}{T} \right)^{1/2} \lambda^{-1} \exp(\lambda x)(1 - \exp(-\lambda T)) \right]^2 \nu(\mathrm{d}u)\,\mathrm{d}x$$
$$+ \int_0^T \int_{\mathbb{R}} u^2 \left[ \left( \frac{2\lambda}{T} \right)^{1/2} \lambda^{-1} \exp(\lambda x)(\exp(-\lambda x) - \exp(-\lambda T)) \right]^2 \nu(\mathrm{d}u)\,\mathrm{d}x$$
$$= \frac{2\lambda^{-1}}{T} \int_{-\infty}^0 [\exp(\lambda x)(1 - \exp(-\lambda T))]^2 \,\mathrm{d}x$$
$$+ \frac{2\lambda^{-1}}{T} \int_0^T [\exp(\lambda x)(\exp(-\lambda x) - \exp(-\lambda T))]^2 \,\mathrm{d}x$$
$$\to 2\lambda^{-1}, \quad \text{as } T \to \infty.$$

According to Theorem 3, to prove (59) it is now sufficient to show that as $T \to \infty$,

$$\int_{-\infty}^T \int_{\mathbb{R}} |u|^3 \left| \left( \frac{2\lambda}{T} \right)^{1/2} \int_{x \vee 0}^T \exp(-\lambda(t-x)) \,\mathrm{d}t \right|^3 \nu(\mathrm{d}u)\,\mathrm{d}x \to 0.$$

But,

$$\int_{-\infty}^{T}\int_{\mathbb{R}}|u|^{3}\left|\left(\frac{2\lambda}{T}\right)^{1/2}\int_{x\vee 0}^{T}\exp(-\lambda(t-x))\,\mathrm{d}t\right|^{3}\nu(\mathrm{d}u)\,\mathrm{d}x$$

$$=\int_{\mathbb{R}}|u|^{3}\nu(\mathrm{d}u)\times\left(\frac{2\lambda^{-1}}{T}\right)^{3/2}\left[\int_{-\infty}^{0}[\exp(\lambda x)(1-\exp(-\lambda T))]^{3}\,\mathrm{d}x \right. \tag{60}$$

$$\left. +\int_{0}^{T}[\exp(\lambda x)(\exp(-\lambda x)-\exp(-\lambda T))]^{3}\,\mathrm{d}x\right]$$

$$\leq \int_{\mathbb{R}}|u|^{3}\nu(\mathrm{d}u)\times\left(\frac{2\lambda^{-1}}{T}\right)^{3/2}\left[\frac{1}{3\lambda}+T\right]\to 0, \tag{61}$$

thus concluding the proof. □

A CLT analogous to (59) also holds in the more general case of an *extended Ornstein–Uhlenbeck Lévy process* $Y_t^{\phi,\lambda}$, defined as

$$Y_t^{\phi,\lambda}=\sqrt{2\lambda}\int_{-\infty}^{t}\int_{\mathbb{R}}u\exp(-\lambda(t-x))\phi(t,x)\widehat{N}(\mathrm{d}u,\mathrm{d}x), \qquad t\geq 0, \tag{62}$$

where $\lambda > 0$ and the kernel $\phi(t,x)$, from $\mathbb{R}_+ \times \mathbb{R}$ to $\mathbb{R}_+$, is such that $0 < \varepsilon < \phi(t,x) < \eta < \infty$ for some strictly positive finite constants $\varepsilon$ and $\eta$. Indeed, in this case, one can easily show that the variance of the random variable $\int_0^T Y_t^{\phi,\lambda}\,\mathrm{d}t$, denoted $V_{\phi,\lambda}(T)$, is an increasing function (whose explicit expression depends on the choice of $\lambda$ and $\phi$) such that

$$a(\varepsilon,\eta,\lambda)\leq\frac{V_{\phi,\lambda}(T)}{T}\leq b(\varepsilon,\eta,\lambda), \qquad T>0, \tag{63}$$

where $a(\varepsilon,\eta,\lambda)$ and $b(\varepsilon,\eta,\lambda)$ are two positive constants not depending on $T$. Finally, the combination of the estimate (63) and of the arguments displayed in formulae (60) and (61) yields the CLT

$$\frac{1}{V_{\phi,\lambda}(T)^{1/2}}\int_{0}^{T}Y_t^{\phi,\lambda}\,\mathrm{d}t\stackrel{\mathrm{law}}{\to}\mathcal{N}(0,1), \qquad T\to\infty.$$

We now turn to quadratic functionals. We shall now suppose that the measure $\nu(\cdot)$ is such that $\int u^j \nu(\mathrm{d}u) < \infty$ for $j=2,4,6$, and $\int u^2 \nu(\mathrm{d}u) = 1$. Plainly, these assumptions yield that

$$\mathbb{E}[(Y_t^{\lambda})^2]=2\lambda\int_{-\infty}^{t}\int_{\mathbb{R}}u^2 \mathrm{e}^{-2\lambda(t-x)}\nu(\mathrm{d}u)\,\mathrm{d}x=1.$$

**Theorem 5.** *For every $\lambda > 0$, define*

$$\hat{h}_t^{\lambda}(u,x;u',x')=2\lambda uu'\exp(-\lambda(t-x)-\lambda(t-x'))\mathbf{1}_{(-\infty,t]^2}(x,x'). \tag{64}$$

*Then as $T \to \infty$:*

1. 
$$\sqrt{T}\left\{\frac{1}{T}\int_0^T I_2(\hat{h}_t^\lambda)\,dt,\ \frac{1}{T}\int_0^T\left[\int_{-\infty}^t\int_{\mathbb{R}} 2\lambda u^2 e^{-2\lambda(t-x)}\widehat{N}(du,dx)\right]dt\right\} \quad (65)$$
$$\overset{\text{law}}{\to} (\mathcal{N}(0,\lambda^{-1}), \mathcal{N}'(0,c_\nu^2)), \quad (66)$$

*where $I_2$ is a double Poisson integral with respect to $\widehat{N}$, and $\mathcal{N}$ and $\mathcal{N}'$ are two independent centered Gaussian random variables with variances given, respectively, by $\lambda^{-1}$ and $c_\nu^2 = \int u^4 \nu(du)$;*

2. 
$$\sqrt{T}\left\{\frac{1}{T}\int_0^T (Y_t^\lambda)^2\,dt - 1\right\} \overset{\text{law}}{\to} \sqrt{\frac{1}{\lambda} + c_\nu^2} \times \mathcal{N}(0,1), \quad (67)$$

*where $\mathcal{N}(0,1)$ indicates a centered standard Gaussian random variable.*

**Proof.** We introduce the notation

$$H_{\lambda,T}(u,x;u',x') = (u \times u')\frac{\mathbf{1}_{(-\infty,T]^2}(x,x')}{T}\{e^{\lambda(x+x')}(1-e^{-2T})\mathbf{1}_{(x \vee x' \le 0)}$$
$$+ e^{\lambda(x+x')}(e^{-2\lambda(x \vee x')} - e^{-2\lambda T})\mathbf{1}_{(x \vee x' > 0)}\}, \quad (68)$$

$$H_{\lambda,T}^*(u,x) = u^2\frac{\mathbf{1}_{(-\infty,T]}(x)}{T}\{e^{2\lambda x}(1-e^{-2T})\mathbf{1}_{(x \le 0)} + e^{2\lambda x}(e^{-2\lambda x} - e^{-2\lambda T})\mathbf{1}_{(x>0)}\}.$$

A standard interchange of deterministic and stochastic integration yields that

$$\sqrt{T}\left\{\frac{1}{T}\int_0^T I_2(\hat{h}_t^\lambda)\,dt,\ \frac{1}{T}\int_0^T\left[\int_{-\infty}^t\int_{\mathbb{R}} 2\lambda u^2 e^{-2\lambda(t-x)}\widehat{N}(du,dx)\right]dt\right\}$$
$$= \{I_2(\sqrt{T}H_{\lambda,T}), I_1(\sqrt{T}H_{\lambda,T})\} \triangleq \{K_2(T), K_1(T)\},$$

where $I_2$ and $I_1$ denote, respectively, a double and a single Wiener–Itô integral with respect to $\widehat{N}$ (observe that, since in (65) the parameter $t$ is integrated with respect to a finite measure, the interchange of deterministic and stochastic integrals can be justified by means of a standard stochastic Fubini theorem – one can, e.g., mimic the proof of [16] Lemma 13, by first approximating the kernels $\hat{h}_t^\lambda$ and $u^2 e^{-2\lambda(t-x)}$ by means of piecewise constant integrands and then by using the isometric properties of single and double Wiener–Itô integrals). Since for every $T$, $K_1(T)$ is a single integral and $K_2(T)$ is a double integral, the joint convergence of the vector $\{K_2(T), K_1(T)\}$ can be studied by means of Theorem 3. By using the same kind of calculations as in the proof of Theorem 4, one easily verifies that

$$\mathbb{E}[K_1(T)^2] = T\int_{\mathbb{R}\times\mathbb{R}} H_{\lambda,T}^*(u,x)^2 \nu(du)\,dx \to c_\nu^2$$

and

$$\int_{\mathbb{R}\times\mathbb{R}} |\sqrt{T}H^*_{\lambda,T}(u,x)|^3 \nu(\mathrm{d}u)\,\mathrm{d}x = T^{3/2} \int_{\mathbb{R}\times\mathbb{R}} |H^*_{\lambda,T}(u,x)|^3 \nu(\mathrm{d}u)\,\mathrm{d}x \sim T^{-1/2} \to 0.$$

We therefore deduce from part A of Theorem 3 that $K_1(T) \stackrel{\text{law}}{\to} N(0,c_\nu^2)$. In view of part B of Theorem 3, the CLT (66) holds if the kernel

$$J_T(u,x;u',x') \triangleq \sqrt{T}H_{\lambda,T}(u,x;u',x') \tag{69}$$

verifies Assumption N (i.e., relations (21)–(23)) and (24). Start by observing that (N-ii) = (22) holds (up to a different normalization) because

$$2\|J_T\|^2_{L^2((\mathbb{R}\times\mathbb{R})^2,(d\nu\times\mathrm{d}x)^2)}$$

$$= 2T\|H_{\lambda,T}\|^2_{L^2((\mathbb{R}\times\mathbb{R})^2,(d\nu\times\mathrm{d}x)^2)}$$

$$= \left(\int_{\mathbb{R}} u^2\nu(\mathrm{d}u)\right)^2 \int_{\mathbb{R}\times\mathbb{R}} \frac{\mathbf{1}_{(-\infty,T]^2}(x,x')}{T^2}\{e^{\lambda(x+x')}(1-e^{-2T})\mathbf{1}_{(x\vee x'\leq 0)}$$

$$\times e^{\lambda(x+x')}(e^{-2\lambda(x\vee x')}-e^{-2\lambda T})\mathbf{1}_{(x\vee x'>0)}\}^2\,\mathrm{d}x\,\mathrm{d}x'$$

$$\longrightarrow \lambda^{-1}.$$

It is therefore sufficient to verify that kernel $J_T$ defined in (69) verifies conditions (N-i) = (21), (N-i) = (23) and the contraction condition (24), namely that

$$\int_{\mathbb{R}\times\mathbb{R}} J_T(\cdot;u,x)^2 \nu(\mathrm{d}u)\,\mathrm{d}x \in L^2(\mathbb{R}\times\mathbb{R},\mathrm{d}\nu\,\mathrm{d}x) \qquad \forall T>0, \tag{70}$$

$$\left\{\int_{\mathbb{R}\times\mathbb{R}} J_T(\cdot;u,x)^4 \nu(\mathrm{d}u)\,\mathrm{d}x\right\}^{1/2} \in L^1(\mathbb{R}\times\mathbb{R},\mathrm{d}\nu\,\mathrm{d}x) \qquad \forall T>0, \tag{71}$$

$$\int_{(\mathbb{R}\times\mathbb{R})^2} J_T(u,x;u',x')^4 \nu(\mathrm{d}u)\,\mathrm{d}x\nu(\mathrm{d}u')\,\mathrm{d}x' \to 0, \qquad T\to\infty, \tag{72}$$

$$\|J_T \star^1_2 J_T\|^2_{L^2(\mathbb{R}\times\mathbb{R},\mathrm{d}\nu\,\mathrm{d}x)} \to 0, \qquad T\to\infty, \tag{73}$$

$$\|J_T \star^1_1 J_T\|^2_{L^2((\mathbb{R}\times\mathbb{R})^2,(\mathrm{d}\nu\,\mathrm{d}x)^2)} \to 0, \qquad T\to\infty. \tag{74}$$

In view of (69), (70) and (71) are implied by the definition of $H_{\lambda,T}$. Relation (72) can be deduced from the relation

$$T^2 \times \|H^2_{\lambda,T}\|^2_{L^2(\mathbb{R}\times\mathbb{R},\mathrm{d}\nu\,\mathrm{d}x)} \sim \frac{1}{T} \to 0. \tag{75}$$

Relations (73) and (74) are a consequence of the two asymptotic relations

$$T^2 \times \int_{\mathbb{R}\times\mathbb{R}} \left(\int_{\mathbb{R}\times\mathbb{R}} H_{\lambda,T}(u,x;u',x')^2 \nu(\mathrm{d}u)\,\mathrm{d}x\right)^2 \nu(\mathrm{d}u')\,\mathrm{d}x'$$

$$\sim \frac{1}{T} \to 0, \tag{76}$$

$$T^2 \int_{(\mathbb{R}\times\mathbb{R})^2} \left( \int_{\mathbb{R}\times\mathbb{R}} H_{\lambda,T}(v,z;u,x) H_{\lambda,T}(v,z;u',x') \nu(\mathrm{d}v) \, \mathrm{d}z \right)^2 \nu(\mathrm{d}u) \, \mathrm{d}x \nu(\mathrm{d}u') \, \mathrm{d}x' \tag{77}$$

$$\sim \frac{1}{T} \to 0.$$

Note that (75)–(77) can be easily checked by resorting to the explicit definition of $H_{\lambda,T}$, as given in (68). The conclusion of point 1 now follows. Point 2 can be deduced from the relation

$$\sqrt{T}\left\{ \frac{1}{T} \int_0^T (Y_t^\lambda)^2 \, \mathrm{d}t - 1 \right\}$$
$$= \sqrt{T}\left\{ \frac{1}{T} \int_0^T I_2(\hat{h}_t^\lambda) \, \mathrm{d}t + \frac{1}{T} \int_0^T \left[ \int_{-\infty}^t \int_\mathbb{R} 2\lambda u^2 \mathrm{e}^{-2\lambda(t-x)} \widehat{N}(\mathrm{d}u, \mathrm{d}x) \right] \mathrm{d}t \right\},$$

which is a consequence of the multiplication formula (14) (in the case $p = q = 1$), applied for every $t$ to the variable $(Y_t^\lambda)^2$. □

*Remark.* In the Poisson case considered here, the double and single integrals in (66) both converge to a non-degenerate Gaussian distribution. This situation is different when the random measure is Gaussian. In that case, the single integral – which corresponds to the contribution of the diagonal – is deterministic.

We can also prove the following CLT for the sample variance of $Y_t^\lambda$.

**Corollary 6.** *With the same notation as in Theorem 5, for every $\lambda > 0$,*

$$\sqrt{T}\left\{ \frac{1}{T} \int_0^T \left( Y_t^\lambda - \frac{1}{T} \int_0^T Y_u^\lambda \, \mathrm{d}u \right)^2 \mathrm{d}t - 1 \right\} \overset{\mathrm{law}}{\to} \sqrt{\frac{1}{\lambda} + c_\nu^2} \times \mathcal{N}(0,1).$$

**Proof.** We simply write $\frac{1}{T} \int_0^T (Y_t^\lambda - \frac{1}{T} \int_0^T Y_u^\lambda \, \mathrm{d}u)^2 \, \mathrm{d}t = \frac{1}{T} \int_0^T (Y_t^\lambda)^2 \, \mathrm{d}t - (\frac{1}{T} \int_0^T Y_t^\lambda \, \mathrm{d}t)^2$ and observe that by (59), $\sqrt{T}(\frac{1}{T} \int_0^T Y_t^\lambda \, \mathrm{d}t)^2 = O_\mathbb{P}(T^{-1/2})$. □

### 4.2. Functionals of random hazard rates in nonparametric Bayesian survival analysis

We shall now discuss some applications to *random hazard rate models* in nonparametric Bayesian survival analysis. These random hazard rates are often represented as generalized Volterra processes of the kind described above. In [5] and [15], the linear and quadratic functionals associated with random hazard rates in some popular Bayesian models are studied by means of the techniques developed in this paper. In what follows, we present the main elements of a Bayesian random hazard rate model, as well as some crucial examples and motivations taken from [5] and [15]. The reader is referred to [5], [10], [15] and the references therein for more details about Bayesian models of survival analysis.

Let $(Z, \mathcal{Z})$ be a measurable space and let $\mu$ be a $\sigma$-finite measure over $(Z, \mathcal{Z})$. As before, we use the notation $\mathcal{Z}_\mu = \{B \in \mathcal{Z} : \mu(B) < \infty\}$. A collection of random variables $N = \{N(B) : B \in \mathcal{Z}_\mu\}$ is called a *non-compensated Poisson measure* with control measure $\mu$ if there exists a compensated Poisson measure $\widehat{N} = \{\widehat{N}(B) : B \in \mathcal{Z}_\mu\}$ (see Section 2) such that

$$N(B) = \widehat{N}(B) + \mu(B) \qquad \forall B \in \mathcal{Z}_\mu. \tag{78}$$

In other words, $N = \{N(B) : B \in \mathcal{Z}_\mu\}$ is a non-compensated Poisson measure if and only if (i) for every $B, C \in \mathcal{Z}_\mu$ such that $B \cap C = \varnothing$, $N(B)$ and $N(C)$ are independent, and (ii) for every $B \in \mathcal{Z}_\mu$,

$$N(B) \stackrel{\text{law}}{=} \mathfrak{P}(B),$$

where $\mathfrak{P}(B)$ is a Poisson random variable with parameter $\mu(B)$.

The basic ingredients of a Bayesian random hazard rate model are the following:

- a *random hazard rate* $\widetilde{h}$, which is a positive generalized Volterra process of the type

$$\widetilde{h}(t) = \int_\mathbb{R} \int_{\mathbb{R}_+} uk(t,x) N(\mathrm{d}u, \mathrm{d}x), \qquad t \geq 0, \tag{79}$$

where $k(t,x) \geq 0$ and $N$ is a suitable non-compensated Poisson measure over $\mathbb{R} \times \mathbb{R}_+$, verifying

$$\lim_{T \to \infty} \int_0^T \widetilde{h}(t) \, \mathrm{d}t = +\infty, \qquad \text{a.s.-}\mathbb{P};$$

- a *random density* with support in $\mathbb{R}_+$, given by

$$f(t) = \widetilde{h}(t) \exp\left\{-\int_0^t \widetilde{h}(s) \, \mathrm{d}s\right\} = \widetilde{h}(t) \exp\{-\widetilde{H}(t)\}, \qquad t \geq 0, \tag{80}$$

where

$$\widetilde{H}(t) = \int_0^t \widetilde{h}(s) \, \mathrm{d}s, \qquad t \geq 0;$$

- a sequence of positive absolutely continuous exchangeable random variables $\mathbf{U} = \{U_n : n \geq 1\}$, representing the *lifetimes* associated with a given population, such that, conditionally on the density $f$ in (80), $\mathbf{U}$ is composed of i.i.d. random variables with common law given by $f$.

Plainly, the initial choice of the law of $\widetilde{h}$ and $f$ (which, in Bayesian terms, is the *prior specification* of the model) is completely encoded by the choices of the control measure of $N$ and of the kernel $k$. Note that (80) gives the following heuristic characterization of $\widetilde{h}(t)$:

$$\widetilde{h}(t)\,dt = \mathbb{P}(t \leq U_1 \leq t + dt | U_1 \geq t, f),$$

meaning that given $f$, the quantity $\widetilde{h}(t)\,dt$ is the probability that the lifetime $U_1$ (or, for that matter, $U_n$) falls in the interval $[t, t+dt]$, conditionally on the fact that $U_1$ is greater than $t$.

Popular choices for $k$ are the following: (a) the *Dykstra–Laud kernel* $k(t,x) = \mathbf{1}_{0 \leq x \leq t}$; (b) the *rectangular kernel* $k(t,x) = \mathbf{1}_{|x-t| \leq \tau}$, where $\tau > 0$ is called the 'bandwith' of the kernel; (c) the Ornstein–Uhlenbeck kernel $k(t,x) = \sqrt{2\lambda}\exp\{-\lambda(t-x)\}\mathbf{1}_{0 \leq x \leq t}$ (note the difference with (58), where the kernel is indeed $\sqrt{2\lambda}\exp\{-\lambda(t-x)\}\mathbf{1}_{-\infty < x \leq t}$).

Popular choices for the control measure of $N$ are:

(i) *generalized Gamma controls* of the type

$$\mu(du, dx) = \Gamma(1-\sigma)^{-1} \exp(-\gamma u) u^{-1-\sigma} \mathbf{1}_{u,x>0}\,du\,dx, \tag{81}$$

where $\sigma \in (0,1)$, $\gamma > 0$ and $\Gamma$ is the usual Gamma function;

(ii) *extended Gamma controls* of the kind

$$\mu(du, dx) = \exp(-\beta(x)u) u^{-1} \mathbf{1}_{u,x>0}\,du\,dx, \tag{82}$$

where $\beta$ is a strictly positive function on $\mathbb{R}$;

(iii) *Beta controls* of the type

$$\mu(du, dx) = (1-u)^{c(x)-1} c(x) \mathbf{1}_{u \in (0,1)} \mathbf{1}_{x>0}\,du\,dx, \tag{83}$$

where $c$ is a strictly positive function on $\mathbb{R}$.

Note that, in general, extended Gamma controls and Beta controls are *non-homogeneous*, in the sense that they cannot be represented in the form $\mu(du, dx) = \nu(du)\,dx$, for some $\sigma$-finite measure $\nu$ on $\mathbb{R}_+$.

The crucial point is that due to the very complex nature of an object such as (79), very little is known about the effect that different parametric choices in (81)–(83) may have on the distributional behavior of $\widetilde{h}$, in particular with respect to functionals of statistical relevance. The idea developed in [15] and, later, in [5], is that one can always represent $\widetilde{h}$ in terms of some underlying *compensated* Poisson measure (via relation (78)) so that Theorem 2 and Theorem 3 of the present paper can be applied in order to obtain CLTs for linear and quadratic functionals of $\widetilde{h}$. The key step is, of course, to represent linear and quadratic functionals as a sum of a single and a double Wiener–Itô integral, a task that can be easily performed by using, for example, the multiplication formula (14). It turns out that in most cases, the constants involved in the CLTs obtained in this way can be expressed very neatly in terms of the different parameters composing the control measure of $N$ and the kernel $k$, thus giving a rough description of the overall 'shape' of the hazard rate $\widetilde{h}$ as a function of the prior specification of the model. As argued in [15] and [5], these kind of results may serve as a guide in the prior analysis since they provide a quite direct way to incorporate prior knowledge into the specification of the law of $\widetilde{h}$.

As an example, we present some results, proved in [15] by means of Theorem 2 and Theorem 3, involving linear and quadratic functionals of the rectangular kernel, under different parametric choices of the control measure associated with $N$. The first statement concerns linear functionals. Recall that the random variables $\widetilde{H}(T) = \int_0^T \widetilde{h}(s)\,\mathrm{d}s$, $T > 0$, have been defined in (80); the symbol $\mathcal{N}(0,c)$ denotes a centered Gaussian random variable with variance $c$.

**Theorem 7 (See [15]).** (1) *Let $N$ have a homogeneous control measure of the type $\mu(\mathrm{d}u, \mathrm{d}x) = \nu(\mathrm{d}u)\,\mathrm{d}x$, where the $\sigma$-finite measure $\nu$ is such that $K_\nu^{(i)} = \int_{\mathbb{R}_+} u^i \nu(\mathrm{d}u) < \infty$, $i = 1, 2$. Let $\widetilde{h}(t)$ be defined via (79), with the rectangular kernel $k(t,x) = \mathbf{1}_{|t-x| \leq \tau}$, for some $\tau > 0$. Then*

$$\frac{1}{\sqrt{T}}[\widetilde{H}(T) - 2\tau K_\nu^{(1)} T] \xrightarrow{\mathrm{law}} \mathcal{N}(0, 4\tau^2 K_\nu^{(2)}).$$

(2) *Let $N$ have a non-homogeneous control measure of the type (82), with $\beta(x) = 1 + x^{1/2}$. Let $\widetilde{h}(t)$ be defined via (79), with the rectangular kernel of bandwith one given by $k(t,x) = \mathbf{1}_{|t-x| \leq 1}$. Then*

$$\frac{1}{\sqrt{\log T}}[\widetilde{H}(T) - 4T^{1/2}] \xrightarrow{\mathrm{law}} \mathcal{N}(0, 4).$$

(3) *Let $N$ have a non-homogeneous control measure of the type (83), with $c(x) \sim x^{1/2}$ as $x \to \infty$. Let $\widetilde{h}(t)$ be defined via (79) with the rectangular kernel $k(t,x) = \mathbf{1}_{|t-x| \leq 1}$. Then*

$$\frac{1}{T^{1/4}}[\widetilde{H}(T) - 2T] \xrightarrow{\mathrm{law}} \mathcal{N}(0, 8).$$

The following statement involves quadratic functionals of rectangular random hazard rates, under a homogeneous assumption on the control measure of $N$. Further results, involving non-homogeneous random measures, are contained at the end of [15] Section 4.2.

**Theorem 8 (See [15], Section 4.2).** *Let $N$ have a homogeneous control measure of the type $\mu(\mathrm{d}u, \mathrm{d}x) = \nu(\mathrm{d}u)\,\mathrm{d}x$, where the $\sigma$-finite measure $\nu$ is such that $K_\nu^{(i)} = \int_{\mathbb{R}_+} u^i \nu(\mathrm{d}u) < \infty$, $i = 1, \ldots, 4$. Let $\widetilde{h}(t)$ be defined via (79), with the rectangular kernel $k(t,x) = \mathbf{1}_{|t-x| \leq \tau}$, for some $\tau > 0$. Then*

$$\sqrt{T}\left[\frac{1}{T}\int_0^T \widetilde{h}(t)^2\,\mathrm{d}t - (2\tau K_\nu^{(2)} + 4\tau^2 (K_\nu^{(1)})^2)\right] \xrightarrow{\mathrm{law}} \mathcal{N}(0, c_1)$$

with $c_1 = 16\tau^2[K_\nu^{(4)}/4 + \tau K_\nu^{(1)} K_\nu^{(3)} + 2\tau[K_\nu^{(2)}]^2/3 + \tau^2[K_\nu^{(2)}]^2 K_\nu^{(1)}]$. *Moreover,*

$$\sqrt{T}\left[\frac{1}{T}\int_0^T \left(\widetilde{h}(t) - \frac{\widetilde{H}(T)}{T}\right)^2 dt - 2\tau K_\nu^{(2)}\right] \xrightarrow{\text{law}} \mathcal{N}(0, c_2),$$

*where* $c_2 = 4\tau^2[K_\nu^{(4)} + 8\tau[K_\nu^{(2)}]^2/3]$.

**Remark.** The paper [5] continues the analysis contained in [15] by establishing *posterior CLTs* for linear and quadratic functionals of random hazard rates, that is, limit theorems involving the law obtained by conditioning on an arbitrary sample of observations $(U_1, \ldots, U_n)$. One of the main findings in [5] is that the posterior CLTs coincide with the prior ones for most models commonly used in Bayesian analysis. These results are then compared with another asymptotic characterization of Bayesian models, known as *consistency* (see, e.g., Drăghici and Ramamoorthi [7]).

# Acknowledgements

This research was supported in part by NSF grants DMS-05-0547 and DMS-07-06786 at Boston University.

# References

[1] Bhansalia, R.J., Giraitis, L. and Kokoszka, P.S. (2007). Approximations and limit theory for quadratic forms of linear processes. *Stochastic Process. Appl.* **117** 71–95. MR2287104
[2] Bhansalia, R.J., Giraitis, L. and Kokoszka, P.S. (2007). Convergence of quadratic forms with nonvanishing diagonals. *Statist. Probab. Lett.* **77** 726–734. MR2356512
[3] Bandorff-Nielsen, O.E. and Shepard, N. (2001). Non Gaussian Ornstein–Uhlenbeck-based models and some of their uses in financial economics. *J. Roy. Statist. Soc. Ser. B* **63** 167–241. MR1841412
[4] Cohen, S. and Taqqu, M.S. (2004). Small and large scale behavior of the Poissonized Telecom process. *Methodol. Comput. Appl. Probab.* **6** 363–379. MR2108557
[5] De Blasi, P., Peccati, G. and Prünster, I. (2007). Asymptotics of posterior hazards. Preprint.
[6] De Jong, P. (1987). A central limit theorem for generalized quadratic forms. *Probab. Theory Related Fields* **75** 261–277. MR0885466
[7] Drăghici, L. and Ramamoorthi, R.V. (2003). Consistency of Dykstra–Laud priors. *Sankhyā* **65** 464–481. MR2028910
[8] Dudley, R.M. (2002). *Real Analysis and Probability.* Cambridge: Cambridge Univ. Press. MR1932358
[9] Hu, Y. and Nualart, D. (2005). Renormalized self-intersection local time for fractional Brownian motion. *Ann. Probab.* **33** 948–983. MR2135309
[10] James, L.F. (2005). Bayesian Poisson process partition calculus with an application to Bayesian Lévy moving averages. *Ann. Statist.* **33** 1771–1799. MR2166562
[11] Kabanov, Y. (1975). On extended stochastic integrals. *Theory Probab. Appl.* **20** 710–722. MR0397877

[12] Marinucci, D. and Peccati, G. (2007). High-frequency asymptotics for subordinated stationary fields on an Abelian compact group. *Stochastic Process. Appl.* **118** 585–613.

[13] Nualart, D. and Peccati, G. (2005). Central limit theorems for sequences of multiple stochastic integrals. *Ann. Probab.* **33** 177–193. MR2118863

[14] Ogura, H. (1972). Orthogonal functionals of the Poisson process. *IEEE Trans. Inform. Theory IT* **18** 473–481. MR0404572

[15] Peccati, G. and Prünster, I. (2006). Linear and quadratic functionals of random hazard rates: An asymptotic analysis. *Ann. Appl. Probab.* To appear.

[16] Peccati, G. (2001). On the convergence of multiple random integrals. *Studia Sci. Math. Hungar.* **37** 429–470. MR1874695

[17] Peccati, G. and Taqqu, M.S. (2007). Stable convergence of $L^2$ generalized stochastic integrals and the principle of conditioning. *Electron. J. Probab.* **12** 447–480. MR2299924

[18] Peccati, G. and Taqqu, M.S. (2006). Stable convergence of multiple Wiener–Itô integrals. *J. Theoret. Probab.* To appear.

[19] Peccati, G. and Taqqu, M.S. (2007). Limit theorems for multiple stochastic integrals. Preprint. Available at http://www.geocities.com/giovannipeccati/. MR2299924

[20] Peccati, G. and Tudor, C.A. (2004). Gaussian limits for vector-valued multiple stochastic integrals. *Séminaire de Probabilités XXXVIII* 247–262. *Lecture Notes in Math.* **1857**. Berlin: Springer. MR2126978

[21] Rota, G.-C. and Wallstrom, C. (1997). Stochastic integrals: A combinatorial approach. *Ann. Probab.* **25**. 1257–1283. MR1457619

[22] Sato, K.-I. (1999). *Lévy Processes and Infinitely Divisible Distributions. Cambridge Studies in Advanced Mathematics* **68**. Cambridge: Cambridge Univ. Press. MR1739520

[23] Surgailis, D. (1984). On multiple Poisson stochastic integrals and associated Markov semigroups. *Probab. Math. Statist.* **3** 217–239. MR0764148

[24] Surgailis, D. (2000). CLTs for Polynomials of linear sequences: Diagram formulae with applications. In *Theory and Applications of Long-Range Dependence* 111–128. Boston: Birkhäuser. MR1956046

[25] Surgailis, D. (2000). Non-CLT's: $U$-statistics, multinomial formula and approximations of multiple Wiener–Itô integrals. In *Theory and Applications of Long-Range Dependence* 129–142. Boston: Birkhäuser. MR1956047

[26] Wolpert, R.L. and Taqqu, M.S. (2005). Fractional Ornstein–Uhlenbeck Lévy processes and the telecom process: Upstairs and downstairs. *Signal Processing* **85** 1523–1545.

[27] Xue, X.-H. (1991). On the principle of conditioning and convergence to mixtures of distributions for sums of dependent random variables. *Stochastic Process. Appl.* **37** 175–186. MR1102868

[28] Yosida, K. (1980). *Functional Analysis*, 6th ed. Berlin: Springer. MR0617913